\documentclass[12pt,a4paper]{amsart}

\usepackage{amssymb, tikz}

\newtheorem{theorem}{Theorem}[section]
\newtheorem{lemma}[theorem]{Lemma}
\newtheorem{prop}[theorem]{Proposition}

\newtheorem*{prop*}{Proposition}
\newtheorem*{lemma*}{Lemma}
\theoremstyle{remark}
\newtheorem{example}[theorem]{Example}
\newtheorem*{example*}{Example}
\newtheorem{remark}[theorem]{Remark}
\newtheorem*{remark*}{Remark}
\newtheorem{rem}[theorem]{Remark}
\newtheorem*{rem*}{Remark}
\newtheorem*{defn*}{Definition}
\newtheorem{defns}[theorem]{Definitions}
\newtheorem{defn}[theorem]{Definition}
\newtheorem{definition}[theorem]{Definition}

\newtheorem*{notation*}{Notation}
\newtheorem{examples}[theorem]{Examples}
\newtheorem{remarks}[theorem]{Remarks}
\newtheorem*{remarks*}{Remarks}

\newtheorem*{thmcomment}{Background}

\newtheorem*{thmcomment2}{Hooptedoodle}

\def\N{{\mathbb N}}\def\Z{{\mathbb Z}}

\newcommand{\addon}{\!\downharpoonleft}

\numberwithin{equation}{section}

\date{\today}

\begin{document}

     \title[Simplicity of $2$-graph algebras associated to Dynamical Systems]{Simplicity of $2$-graph algebras associated to Dynamical Systems}
     \author[Peter lewin and David Pask]{Peter Lewin and David Pask}
     \address{School of Mathematics and Applied Statistics \\
        Austin Keane Building (15) \\
         University of Wollongong \\
         NSW 2522 \\ AUSTRALIA}
         \email{\{dpask,pkl621\}@uow.edu.au}

    \begin{abstract}
We give a combinatorial description of a family of $2$-graphs which subsumes those described by Pask, Raeburn and Weaver.
Each 2-graph $\Lambda$ we consider has an associated $C^*$-algebra, denoted $C^*(\Lambda)$, which is simple and purely infinite when $\Lambda$ is aperiodic. We give new, straightforward conditions which ensure that $\Lambda$ is aperiodic.
These conditions are highly tractable as we only need to consider the finite set of vertices of $\Lambda$ in order to identify aperiodicity. In addition, the path space of each 2-graph can
be realised as a two-dimensional dynamical system which we show must have zero entropy.

    \end{abstract}


    \thanks{This research was supported by ARC Discovery Project DP0665131}


    \maketitle

    \vskip4em

    \section{Introduction.}

    Higher-rank graphs, or $k$-graphs, were introduced by Kumjian and Pask~\cite{kp} as a graph based model for the higher-rank Cuntz-Krieger algebras studied by Robertson and Steger~\cite{robst}. The $C^*$-algebras associated to higher-rank graphs have been generating a lot of interest~\cite{gevans,fmy,shotwell, sims2}. In this paper we show how to obtain a class of $2$-graphs from a small set of parameters, which we will call ``basic data''. Examples
    generated in this way include all those considered by Pask, Raeburn and Weaver in \cite{prw}. Our examples
    will therefore add to the emerging connections between $2$-graphs, classification theory and two-dimensional shift spaces (\cite{kp2,prrs,spiel,sz}).

    This paper has three aims: The first aim is to model a family of $2$-graphs which subsumes those studied in \cite{prw}. The second aim is to identify when the associated $C^*$-algebra is simple and purely infinite using the new results of~\cite{ls, robsi1}. The third aim is to establish connections between the 2-graphs we generate and two-dimensional shift spaces, as defined in \cite{bfm, ms, qt}. In doing this, we generalise  certain results of Pask, Raeburn and Weaver in \cite{prw}.

    For our first aim, we write down a set of basic data, $BD$ which determines a finite set of labelled tiles. In Theorem~\ref{2graph} we show that $BD$ completely specifies a $2$-graph, denoted $\Lambda_{BD}$: The vertices of $\Lambda_{BD}$ are the set of labelled tiles, with an edge
    between two vertices if their labels overlap in a particular way. The $2$-graph $\Lambda_{BD}$ has a visual representation, called the skeleton of $\Lambda_{BD}$. We prove several basic properties about the skeleton of $\Lambda_{BD}$ in Proposition~\ref{sk} and Proposition~\ref{sc}. In particular we show that $\Lambda_{BD}$ is strongly connected, row-finite and has no sources.

    For our second aim, we investigate when the $C^*$-algebra associated to $\Lambda_{BD}$ is purely infinite and simple.
    The results of~\cite{kp,ls, rsy1, shotwell} lead to a characterisation of simplicity of $k$-graph $C^*$-algebras in terms of two conditions on the underlying $k$-graph; namely aperiodicity and cofinality. Since all our examples are cofinal, we focus our attention on identifying when $\Lambda_{BD}$ is aperiodic.
    Theorem~\ref{APT} gives a condition on the set of vertices of $\Lambda_{BD}$ which ensures that $\Lambda_{BD}$ is aperiodic. As the set of vertices is finite, checking this condition is straightforward. Our aperiodicity result then allows us to prove the primary result of this paper, Theorem~\ref{class} which gives necessary conditions for $C^* ( \Lambda_{BD} )$ to be simple and purely infinite.

    For our third aim, we note that the shift map $\sigma$ induces a $\mathbb{Z}^2$-action on the two-sided infinite path space $\Lambda_{BD}^\Delta$ of $\Lambda_{BD}$. In Theorem~\ref{hmeo} we show that $(\Lambda_{BD}^\Delta, \sigma)$ is a two-dimensional shift of finite type,  as studied by Quas and Trow~\cite{qt}. Unlike the $2$-graphs studied in \cite{prw} the path space $\Lambda_{BD}^\Delta$ has no underlying algebraic properties.

    We first introduce some background in Section~\ref{bg}. In Section~\ref{gss} we describe the basic data from which we generate a 2-graph  $\Lambda_{BD}$. We go on to prove some general properties about the skeletons of $\Lambda_{BD}$. In Section~\ref{aps} we detail a condition on $\Lambda_{BD}^0$ which will ensure that
    $\Lambda_{BD}$ is aperiodic. We use this condition to establish when $C^*(\Lambda_{BD} )$ is unital, nuclear, purely infinite and simple. This result is proved in Section~\ref{simplicity}. Finally we connect our work to the literature on two-dimensional shift spaces in Section~\ref{ss}.

    \subsection*{Acknowledgements}

    The first author would like to point out a significant portion of Section~\ref{gss} has been based on the ideas of \cite[Section~3]{prw} and would like to acknowledge Pask, Raeburn and Weaver for providing a preprint of their paper as motivation for this publication. The second author would like to thank Anthony Quas for useful discussions regarding the entropy of shift spaces arising from higher-rank graphs.

    \section{Background Information.}\label{bg}

    \subsection{Higher-rank graphs.}

        We regard $\N^k$ as a semigroup under addition, with basis elements $e_i$. For $m, n \in \N^k$ we write $m_i, n_i$ for the $i$th coordinates of $m$ and $n$, $m \vee n$ for their coordinatewise maximum and $m \wedge n$ for their coordinatewise minimum. We let $\le$ denote the usual coordinatewise partial order on $\N^k$.

    A $k$-graph is a pair $(\Lambda, d)$ consisting of a countable category $\Lambda$, together with a functor $d : \Lambda \rightarrow \N^k$ which satisfies the \emph{factorisation property}: for every $\lambda \in \Lambda$ and $m, n \in \N^k$ satisfying $d(\lambda) = m + n$ there exists unique $\mu, \nu \in \Lambda$ such that $d(\mu) = m$, $d(\nu) = n$ and $\lambda = \mu \nu$. We refer to $d$ as the \emph{degree map} and think of it at a generalised length function for elements of $\Lambda$. We think of objects in the category as vertices, and morphisms as paths. We refer to morphisms of degree $e_i$ as edges. For $n \in \N^k$ let $\Lambda^n := \{ \lambda \in \Lambda : d(\lambda) = n \}$ be the set of all paths of degree $n$. We frequently omit the functor $d$ when referring to a $k$-graph.

    For each $\lambda \in \Lambda^n$ the factorisation property implies there exist unique elements $r(\lambda), s(\lambda) \in \Lambda^0$ such that $r(\lambda) \lambda = \lambda = \lambda s(\lambda)$. Hence we can identify the objects in a $k$-graph with the paths of degree 0. For $X \subseteq \Lambda$ and $v \in \Lambda^0$ let $v X = \{ \lambda \in X : r ( \lambda ) = v \}$ and $X v = \{ \lambda \in  X : s ( \lambda ) = v \}$. For $u, v \in \Lambda^0$ set $u X v = uX \cap Xv$.

        In this paper we will be working exclusively with the case when $k = 2$. All 2-graphs in this paper are row-finite and free of sinks and sources: for every $v \in \Lambda^0$ and every $n \in \N^2$ we have $0 < | v \Lambda^{n}| < \infty$ and $0 < | \Lambda^n v |$.

        To visualise a 2-graph we draw its \emph{skeleton}. The skeleton of a 2-graph $\Lambda$ is a bicoloured directed graph with vertices $\Lambda^0$ and edges $\Lambda^{e_1} \cup \Lambda^{e_2}$ with range and source maps inherited from $\Lambda$. We call edges of degree $e_1$ blue edges, and draw them using solid lines. We call edges of degree $e_2$ red edges, and draw them using dashed lines. The following picture is an example of a skeleton of a 2-graph $\Lambda$.


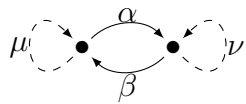
    \[
        \begin{tikzpicture}[scale=1.2]
        \node[circle,inner sep=0.5pt] (00) at (0, 0) {$\bullet$}
                edge[-latex, dashed, loop, out=-135, in=135, min distance=25] node[inner sep=0pt, pos=0.5] (m) {} (00);
        \node[circle,inner sep=0.5pt] (10) at (1, 0) {$\bullet$}
                edge[-latex, dashed, loop, out=-45, in=45, min distance=25] node[inner sep=0pt, pos=0.5] (n) {} (10)
                edge[latex-, out=135, in=45] node[inner sep=0pt, pos=0.5] (a) {} (00)
                edge[-latex, out=-135, in=-45] node[inner sep=0pt, pos=0.5] (b) {} (00);
        \node[inner sep=0.5pt, anchor=south] at (a) {$\alpha$};
        \node[inner sep=0.5pt, anchor=north] at (b) {$\beta$};
        \node[inner sep=0.5pt, anchor=east] at (m) {$\mu$};
        \node[inner sep=0.5pt, anchor=west] at (n) {$\nu$};
        \end{tikzpicture}
        \]

        \noindent The skeleton itself does not encode all the information in $\Lambda$, we also need to record the factorisation of paths of degree $e_1 + e_2$ in $\Lambda$. In the skeleton given above, there
        are only two such factorisations, namely $\alpha \mu = \nu \alpha$ and $\beta \mu = \nu \beta$. Since there is at most one edge of each colour between any two vertices, and for every blue-red path between two vertices there is exactly one red-blue path between the same vertices, the skeleton itself determines the factorisations of paths of degree $e_1+e_2$, and so we need not record them. The 2-graphs
        we study in this paper will also have this property (see Proposition~\ref{sk}). A concise overview of skeletons can be found in \cite[Section~2]{rsy1}.

     Let $\Lambda$ be a row-finite $2$-graph with no sources. A \emph{Cuntz-Krieger $\Lambda$-family} is a set of partial isometries $\{ S_\lambda : \lambda \in \Lambda \}$ satisfying:

        \begin{enumerate}

        \item[(CK1)] $\{ S_v : v \in \Lambda^0 \}$ is a collection of mutually orthogonal projections;

        \item[(CK2)] $S_\mu S_\nu = S_{\mu \nu}$ whenever $s(\mu) = r(\nu)$;

        \item[(CK3)] $S_\lambda^* S_\lambda = S_{s(\lambda)}$ for every $\lambda \in \Lambda$ and;

        \item[(CK4)] $S_v = \sum_{\lambda \in v \Lambda^n} S_\lambda S_\lambda^*$ for every $v \in \Lambda^0$ and every $n \in \N^2$.

        \end{enumerate}

        The $C^*$-algebra $C^* ( \Lambda )$ is then the universal $C^*$-algebra generated by a Cuntz-Krieger $\Lambda$-family $\{ s_\lambda : \lambda \in \Lambda \}$: Given any other Cuntz-Krieger $\Lambda$-family $\{ t_\lambda : \lambda \in \Lambda \}$ in a $C^*$-algebra $B$ there exists a unique homomorphism $\pi_S : C^* ( \Lambda ) \to B$  such that $\pi_S (s_\lambda) = t_\lambda$ for every $\lambda \in \Lambda$. By \cite[Proposition~2.11]{kp} there is a Cuntz-Krieger $\Lambda$-family where each $S_\lambda$ is non-zero. Further background on higher-rank graph $C^*$-algebras can be found in \cite[Section~3]{rsy1}.

    \subsection{Two-dimensional shift spaces.}

    Let $A$ be a finite set of symbols and denote by $A^{\Z^2}$ the set of all functions $\Z^2 \rightarrow A$. We think of elements in $A^{\Z^2}$ as being an infinite two-dimensional array of symbols from $A$. For $x \in A^{\Z^2}$ and $S \subset \Z^2$, let $x|_S$ denote the restriction of $x$ to $S$. A map $E : S \rightarrow A$ is called a \emph{configuration}. For a subset $X$ of $A^{\Z^2}$ a configuration $E$ is said to be \emph{allowed} for $X$ if there exists $x \in X$ such that $x|_S = E$. In this case we say that $E$ occurs in $x$. 

    For every $b \in \Z^2$, define $\sigma_b : A^{\Z^2} \rightarrow A^{\Z^2}$ by $\sigma_b (x) (m) = x(m + b)$ for all $x \in A^{\Z^2}$ and $m \in \Z^2$. Since $\sigma_b$ is a homeomorphism and $\sigma_a \sigma_b = \sigma_{a+b}$ for all $a,b \in \Z^2$ we see that $\Z^2$ acts on $A^{\Z^2}$. A nonempty subset $X$ of $A^{\Z^2}$ which is invariant under $\sigma_b$ for all $b \in \Z^2$ is called a \emph{shift space}. A shift space $X$ is called a \emph{shift of finite type} if there exists a finite set $S \subset \Z^2$ and a nonempty subset $Q \subset A^S$ such that
    \begin{align*}
            X = \{ x \in A^{\Z^2} : x|_{S + b} \in Q ~ \text{for every} ~ b \in \Z^2 \}.
    \end{align*}
We think of $Q$ as a finite set of allowed configurations for $X$.

        For a shift space $X$ and a fixed $d \in \N$, let $\mathcal{B}_d(X)$ denote the set of all  $d \times d$ square configurations which occur in $X$. We define the \emph{topological entropy} of a shift space $X$ to be $\lim_{d \rightarrow \infty} \frac{1}{2^d} \log |\mathcal{B}_d(X)|$ and denote it by $h(X)$. Further background information about the shift spaces used in this paper, and their entropy can be found in~\cite{qt}.

        \section{The 2-graph construction.}\label{gss}

    \noindent We begin by introducing some notation for the data which we use to determine a $2$-graph.

     A subset $T \subset \N^2$ is \emph{hereditary} if $n \in T$ and $0 \leq m \leq n$ implies $m \in T$. Define a \emph{tile} $T$ to be a hereditary subset of $\N^2$ with finite cardinality  $| T |$ and let $(c_1, c_2) := \bigvee \{ i : i \in T \}$. Then, using the visualisation convention of \cite[Section 1]{prw} for tiles, $T$ consists of $|T|$ squares with $c_1+1$ squares in the first row and $c_2+1$ squares in the first column.
If $T$ is a tile with $|T| \ge 2$ then we say that $T$ is \emph{nondegenerate}.

      \begin{example}\label{ex}

    The set $T = \{ 0, e_1, 2e_1, e_1 + e_2, e_2 \}$ is a nondegenerate tile with $c_1=2$ and $c_2=1$. We visualise $T$ as $\begin{tikzpicture}[scale=0.3]
    \draw (1, 0) -- (1, 2) -- (0, 2) -- (0, 0) -- (3, 0) -- (3, 1) -- (0, 1);
    \draw (2, 0) -- (2, 1);
    \draw (1, 2) -- (2, 2) -- (2, 1);
    \end{tikzpicture}$. Note that the first row has  $c_1 + 1=3$ squares and the first column has $c_2 + 1=2$ squares.

    \end{example}

     If $T$ is a nondegenerate tile we define $P = T \backslash \{ c_1 e_1 , c_2 e_2 \}$. If $c_1 , c_2  \ge 1$ then we visualise $P$ as $T$ with the top left and bottom right squares removed. If $c_1 = 0$, $c_2 \ge 1$  (resp.\ $c_2 = 0$, $c_1 \ge 1$) then we visualise $P$ as $T$ with the squares corresponding to the origin $0$ and $c_2 e_2$ (resp.\ $c_1 e_1$) removed. Note that $P = \emptyset$ if and only if $|T|=2$.

    \begin{defn} \label{def:bdef} A set of \emph{basic data}, is a triple $BD := (T, A, \{ f_p : p \in A^{P} \} )$ where $T$ is a nondegenerate tile, $A$ is a nonempty set (which we call the alphabet), and for each $p \in A^P$ the map $f_p : A \rightarrow A$ is a bijection.
    \end{defn}

    \begin{remarks}
    \begin{enumerate}
    \item Recall that for a nonempty set $A$ we have $A^\emptyset = \{ \emptyset \}$.
    \item If $T$ is degenerate, that is $T=\{0\}$, then we
    define the basic data to be $(T,A, \{ f_a \})$ where for $a \in A$, $f_a$ is the function $f_a : T \to A$ given by $f_a (0)=a$. Following this
    case through the following text would obfuscate the main results, we comment separately on this case in Remarks \ref{rem:T01}, \ref{rem:T02} and \ref{rem:T03}.
\end{enumerate}
\end{remarks}

    \noindent The main result of this section is Theorem~\ref{2graph} which describes the unique 2-graph $\Lambda_{BD}$ associated to the basic data $BD$. In order to state this theorem we require some notation. For each $p \in A^P$ and $a \in A$ define $F_{p, a} : T \rightarrow A$ by
    \begin{align} \label{eq:Fpdef}
            F_{p, a} (t) &=
                    \begin{cases}
                            a ~ &\text{when} ~ t = c_2 e_2, \\
                            f_p (a) ~ &\text{when} ~ t = c_1 e_1, \\
                            p(t) ~ &\text{when} ~ t \in P . \\
                    \end{cases}
    \end{align}

\noindent If $P=\emptyset$ then either  $T = \{ 0 , e_1 \}$ in which case $F_{\emptyset,a} (0) = a$, $F_{\emptyset,a} (e_1) = f_\emptyset (a)$ for all $a \in A$ or $T=\{ 0 , e_2 \}$ in which case $F_{\emptyset,a} (e_2) = a$, $F_{\emptyset,a} (0) = f_\emptyset (a)$ for all $a \in A$.

We denote the collection $\bigcup_{p \in A^P, a \in A} F_{p, a}$ by $\Lambda_{BD}^0$.

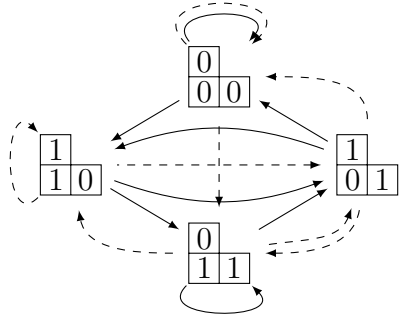
\begin{defns}
    For $S \subset \N^2$ and $n \in \N^2$ let $S + n = \{ m + n : m \in S \}$ denote the translate of $S$ by $n$ and let $T(n) := \bigcup_{0 \leq m \leq n} T + m$ denote the union of translates of $T$ by $m \le n$.
\end{defns}


    \begin{example} For the tile $T = \{ 0, e_1, 2e_1, e_1 + e_2, e_2 \}$ defined in Example~\ref{ex} we have
    \[
    \begin{tikzpicture}[scale=0.3]
\node at (-4,1.3) {$T(2,1)=$};
            \draw (0, 3) -- (0, 0) -- (5, 0) -- (5, 2) -- (4, 2) -- (4, 3) -- (0, 3);
            \draw (0, 1) -- (5, 1);
            \draw (0, 2) -- (4, 2) -- (4, 0);
            \draw (1, 3) -- (1, 0);
            \draw (2, 3) -- (2, 0);
            \draw (3, 3) -- (3, 0);
    \end{tikzpicture}.
    \]
    \end{example}

\noindent If $f : S \rightarrow A$ is a function defined on a subset $S \subset \N^2$ containing $T + n$, then we define\footnote{In \cite{prw} $f \addon_{T+n}$ was denoted $f |_{T+n}$} $f\addon_{T + n} : T \rightarrow A$ by $f \addon_{T + n}(i) = f(i + n) ~ \text{for} ~ i \in T$. Note that $f\addon_T = f|_T$.

 \begin{defns} For $n \in \N^2$ a \emph{path of degree $n$} is a function $\lambda : T(n) \rightarrow A$ such that $\lambda \addon_{T + m} \in \Lambda_{BD}^0$ for $0 \leq m \leq n$, with source $s(\lambda) = \lambda\addon_T$ and range $r ( \lambda ) = \lambda\addon_{T+n}$. We denote the set of all paths of degree $n$ by $\Lambda_{BD}^n$ and define $\Lambda_{BD}^* = \bigcup_{n \in \N^2} \Lambda_{BD}^n$.
 \end{defns}

 \noindent For $\lambda \in \Lambda_{BD}^\ell$ and $0 \leq m \leq n \leq \ell$ the \emph{factorisation} $\lambda (m, n)$ is a path of degree $n - m$ defined by
        \[
                \lambda (m, n) (i) = \lambda (m + i) ~ \text{for} ~ i \in T(n - m).
        \]
Observe in particular that $\lambda (m, m) = \lambda\addon_{T + m} \in \Lambda_{BD}^0$. This is a deviation from notation used elsewhere in the literature \cite{fmy, ls, rsy1, robsi1} and so it is worth highlighting the fact that $\lambda (i) \neq \lambda(i, i)$. This is because $\lambda(i) \in A$ is a value of $\lambda$ at a specific coordinate while $\lambda(i, i) \in A^T$ is a vertex. Thus we continue the convention of \cite{prw}.

        Now we  have all the notation required to state the main theorem of this section.

        \begin{theorem}\label{2graph}

      Given basic data $BD$, the quadruple $\Lambda_{BD} := ( \Lambda_{BD}^0, \Lambda_{BD}^*, r, s)$ is a category. The function $d : \Lambda_{BD} \rightarrow \N^2$ given by $d(\lambda) = n$ for $\lambda \in \Lambda_{BD}^n$ is a functor satisfying the factorisation property; hence $( \Lambda_{BD} , d)$ is a 2-graph.

    \end{theorem}

    \noindent The following two key examples were used to motivate Theorem~\ref{2graph}.

        \begin{examples}\label{ldr}
        \begin{enumerate}
        \item
Let $A = \{ 0, 1 \}$ and $T = \{ 0, e_1, e_2 \}$ be a nondegenerate tile. Then $P=\{0\}$ and so $A^P = \{ p[0] , p[1] \}$ where for $i\in A$ we define $p[i]: P \to A$ by $p[i](0)=i$. To complete our basic data $BD$ we define bijections $f_{p[0]} , f_{p[1]} : A \rightarrow A$ by $f_{p[0]} (a) = a$ and $f_{p[1]} (a) = a + 1$ (mod 2) for all $a \in A$. Then by Theorem~\ref{2graph} the $2$-graph $\Lambda_{BD}$ associated to $BD$ has vertices
    \[
            F_{p[0], 0} =
            \begin{tikzpicture}[scale=0.4]
                    \draw (1, 0) -- (1, 2) -- (0, 2) -- (0, 0) -- (2, 0) -- (2, 1) -- (0, 1);
                    \node[circle,inner sep=0.5pt] () at (0.5, 0.5) {$0$};
                    \node[circle,inner sep=0.5pt] () at (1.5, 0.5) {$0$};
                    \node[circle,inner sep=0.5pt] () at (0.5, 1.5) {$0$};
                    \end{tikzpicture},
              F_{p[0], 1} =
            \begin{tikzpicture}[scale=0.4]
                    \draw (1, 0) -- (1, 2) -- (0, 2) -- (0, 0) -- (2, 0) -- (2, 1) -- (0, 1);
                    \node[circle,inner sep=0.5pt] () at (0.5, 0.5) {$0$};
                    \node[circle,inner sep=0.5pt] () at (1.5, 0.5) {$1$};
                    \node[circle,inner sep=0.5pt] () at (0.5, 1.5) {$1$};
            \end{tikzpicture},
            F_{p[1], 0} =
            \begin{tikzpicture}[scale=0.4]
                    \draw (1, 0) -- (1, 2) -- (0, 2) -- (0, 0) -- (2, 0) -- (2, 1) -- (0, 1);
                    \node[circle,inner sep=0.5pt] () at (0.5, 0.5) {$1$};
                    \node[circle,inner sep=0.5pt] () at (1.5, 0.5) {$1$};
                    \node[circle,inner sep=0.5pt] () at (0.5, 1.5) {$0$};
            \end{tikzpicture},
            F_{p[1], 1} =
            \begin{tikzpicture}[scale=0.4]
                    \draw (1, 0) -- (1, 2) -- (0, 2) -- (0, 0) -- (2, 0) -- (2, 1) -- (0, 1);
                    \node[circle,inner sep=0.5pt] () at (0.5, 0.5) {$1$};
                    \node[circle,inner sep=0.5pt] () at (1.5, 0.5) {$0$};
                    \node[circle,inner sep=0.5pt] () at (0.5, 1.5) {$1$};
            \end{tikzpicture}.
    \]
    The above blocks form the Ledrappier system, named after the two-dimensional shift space studied in~\cite{ledrap}. In~\cite{prw} it was shown that the associated 2-graph has skeleton:
    \[
            \begin{tikzpicture}[scale=0.39]
                    \node[circle,inner sep=0.5pt] (p11) at (1, 4) {\begin{tikzpicture}[scale=0.4]
                                                                                                                \draw (1, 3) -- (1, 5) -- (0, 5) -- (0, 3) -- (2, 3) -- (2, 4) -- (0, 4);
                                                                                                                \end{tikzpicture}}
                            edge[-latex, loop, out=130, in=50, min distance=25, looseness=2] (p11)
                            edge[-latex, loop, out=135, in=45, min distance=35, dashed, looseness=3] (p11)
                            ;
                            \node at (0.5, 3.5) {$0$};
                            \node at (1.5, 3.5) {$0$};
                            \node at (0.5, 4.5) {$0$};
                    \node[circle,inner sep=0.5pt] (p00) at (-4, 1) {\begin{tikzpicture}[scale=0.4]
                                                                                                                \draw (-4, 0) -- (-4, 2) -- (-5, 2) -- (-5, 0) -- (-3, 0) -- (-3, 1) -- (-5, 1);
                                                                                                                \end{tikzpicture}}
                            edge[latex-] (p11)
                            edge[-latex, loop, out=-135, in=135, min distance=35, dashed, looseness=2] (p00);
                            \node at (-4.5, 0.5) {$1$};
                            \node at (-4.5, 1.5) {$1$};
                            \node at (-3.5, 0.5) {$0$};
                    \node[circle,inner sep=0.5pt] (p01) at (1, -2) {\begin{tikzpicture}[scale=0.4]
                                                                                                                \draw (1, -3) -- (1, -1) -- (0, -1) -- (0, -3) -- (2, -3) -- (2, -2) -- (0, -2);
                                                                                                                \end{tikzpicture}}
                            edge[latex-] (p00)
                            edge[-latex, loop, out=-135, in=-45, min distance=25, looseness=2] (p01)
                            edge[-latex, dashed, out=180, in=-80] (p00)
                            edge[latex-, dashed] (p11);
                            \node at (0.5, -2.5) {$1$};
                            \node at (0.5, -1.5) {$0$};
                            \node at (1.5, -2.5) {$1$};
                    \node[circle,inner sep=0.5pt] (p10) at (6, 1) {\begin{tikzpicture}[scale=0.4]
                                                                                                                \draw (6, 0) -- (6, 2) -- (5, 2) -- (5, 0) -- (7, 0) -- (7, 1) -- (5, 1);
                                                                                                                \end{tikzpicture}}
                            edge[-latex] (p11)
                            edge[-latex, out=160, in=20] (p00)
                            edge[latex-, out=-160, in=-20] (p00)
                            edge[latex-] (p01)
                            edge[-latex, dashed, out=90, in=0] (p11)
                            edge[latex-, dashed] (p00)
                            edge[-latex, dashed, out=-100, in=0] (p01)
                            edge[latex-, dashed, out=-110, in=10] (p01);
                            \node at (5.5, 0.5) {$0$};
                            \node at (5.5, 1.5) {$1$};
                            \node at (6.5, 0.5) {$1$};
            \end{tikzpicture}
    \]

   \item Let $T = \{ 0, e_1, e_2, e_1 + e_2 \}$ be a nondegenerate tile and $A = \{ 0, 1 \}$. Then $P=\{0,e_1+e_2\}$ and so $A^P = \{ p[0,0], p[0,1], p[1,0],p[1,1] \}$ where for $a, b \in A$ we define $p[a, b] : P \rightarrow A$ by $p[a, b](0) = a$ and $p[a, b](e_1 + e_2) = b$.
To obtain a set of basic data $BD$ we define a set of bijections as follows: let $f_{p[0, 0]} (a) = a$ for every $a \in A$ and $f_{p} (a) = a + 1$ (mod 2) for every $a \in A$ and every $p \in A^P \setminus \{ p[0, 0] \}$. Then by Theorem~\ref{2graph} the $2$-graph $\Lambda_{BD}$ associated to this basic data has vertices
\begin{equation} \label{eq:8tiles}
        \begin{array}{cccc}
        F_{p[0, 0], 0} = \begin{tikzpicture}[scale=0.4]
                \draw (0, 0) -- (2, 0) -- (2, 2) -- (0, 2) -- (0, 0);
                \draw (1, 0) -- (1, 2);
                \draw (0, 1) -- (2, 1);
                \node[circle,inner sep=0.5pt] () at (0.5, 0.5) {$0$};
                \node[circle,inner sep=0.5pt] () at (1.5, 0.5) {$0$};
                \node[circle,inner sep=0.5pt] () at (0.5, 1.5) {$0$};
                \node[circle,inner sep=0.5pt] () at (1.5, 1.5) {$0$};
    \end{tikzpicture} , &
    F_{p[0, 0], 1} = \begin{tikzpicture}[scale=0.4]
                \draw (0, 0) -- (2, 0) -- (2, 2) -- (0, 2) -- (0, 0);
                \draw (1, 0) -- (1, 2);
                \draw (0, 1) -- (2, 1);
                \node[circle,inner sep=0.5pt] () at (0.5, 0.5) {$0$};
                \node[circle,inner sep=0.5pt] () at (1.5, 0.5) {$1$};
                \node[circle,inner sep=0.5pt] () at (0.5, 1.5) {$1$};
                \node[circle,inner sep=0.5pt] () at (1.5, 1.5) {$0$};
    \end{tikzpicture} , &
    F_{p[0, 1], 0} = \begin{tikzpicture}[scale=0.4]
                \draw (0, 0) -- (2, 0) -- (2, 2) -- (0, 2) -- (0, 0);
                \draw (1, 0) -- (1, 2);
                \draw (0, 1) -- (2, 1);
                \node[circle,inner sep=0.5pt] () at (0.5, 0.5) {$0$};
                \node[circle,inner sep=0.5pt] () at (1.5, 0.5) {$0$};
                \node[circle,inner sep=0.5pt] () at (0.5, 1.5) {$1$};
                \node[circle,inner sep=0.5pt] () at (1.5, 1.5) {$1$};
    \end{tikzpicture} , &
    F_{p[0, 1], 1} = \begin{tikzpicture}[scale=0.4]
                \draw (0, 0) -- (2, 0) -- (2, 2) -- (0, 2) -- (0, 0);
                \draw (1, 0) -- (1, 2);
                \draw (0, 1) -- (2, 1);
                \node[circle,inner sep=0.5pt] () at (0.5, 0.5) {$0$};
                \node[circle,inner sep=0.5pt] () at (1.5, 0.5) {$1$};
                \node[circle,inner sep=0.5pt] () at (0.5, 1.5) {$0$};
                \node[circle,inner sep=0.5pt] () at (1.5, 1.5) {$1$};
    \end{tikzpicture} , \\
    F_{p[1, 0], 0} = \begin{tikzpicture}[scale=0.4]
                \draw (0, 0) -- (2, 0) -- (2, 2) -- (0, 2) -- (0, 0);
                \draw (1, 0) -- (1, 2);
                \draw (0, 1) -- (2, 1);
                \node[circle,inner sep=0.5pt] () at (0.5, 0.5) {$1$};
                \node[circle,inner sep=0.5pt] () at (1.5, 0.5) {$0$};
                \node[circle,inner sep=0.5pt] () at (0.5, 1.5) {$1$};
                \node[circle,inner sep=0.5pt] () at (1.5, 1.5) {$0$};
    \end{tikzpicture} , &
    F_{p[1, 0], 1} = \begin{tikzpicture}[scale=0.4]
                \draw (0, 0) -- (2, 0) -- (2, 2) -- (0, 2) -- (0, 0);
                \draw (1, 0) -- (1, 2);
                \draw (0, 1) -- (2, 1);
                \node[circle,inner sep=0.5pt] () at (0.5, 0.5) {$1$};
                \node[circle,inner sep=0.5pt] () at (1.5, 0.5) {$1$};
                \node[circle,inner sep=0.5pt] () at (0.5, 1.5) {$0$};
                \node[circle,inner sep=0.5pt] () at (1.5, 1.5) {$0$};
    \end{tikzpicture} , &
    F_{p[1, 1], 0} = \begin{tikzpicture}[scale=0.4]
                \draw (0, 0) -- (2, 0) -- (2, 2) -- (0, 2) -- (0, 0);
                \draw (1, 0) -- (1, 2);
                \draw (0, 1) -- (2, 1);
                \node[circle,inner sep=0.5pt] () at (0.5, 0.5) {$1$};
                \node[circle,inner sep=0.5pt] () at (1.5, 0.5) {$0$};
                \node[circle,inner sep=0.5pt] () at (0.5, 1.5) {$1$};
                \node[circle,inner sep=0.5pt] () at (1.5, 1.5) {$1$};
    \end{tikzpicture} , &
    F_{p[1, 1], 1} = \begin{tikzpicture}[scale=0.4]
                \draw (0, 0) -- (2, 0) -- (2, 2) -- (0, 2) -- (0, 0);
                \draw (1, 0) -- (1, 2);
                \draw (0, 1) -- (2, 1);
                \node[circle,inner sep=0.5pt] () at (0.5, 0.5) {$1$};
                \node[circle,inner sep=0.5pt] () at (1.5, 0.5) {$1$};
                \node[circle,inner sep=0.5pt] () at (0.5, 1.5) {$0$};
                \node[circle,inner sep=0.5pt] () at (1.5, 1.5) {$1$};
    \end{tikzpicture}.
        \end{array}
        \end{equation}
\end{enumerate}
        \end{examples}


\begin{rem}
        In \cite[Theorem 3.4]{prw} Pask, Raeburn and Weaver use four parameters to generate a 2-graph: a tile $T$, an alphabet $\Z / q \Z$ where $q \geq 2$, a trace $t \in \Z / q \Z$ and a rule $w : T \to \Z / q \Z$ such that $w ( c_1 e_1 )$ and $w ( c_2 e_2 )$ are invertible in the ring $\Z / q \Z$. The vertices of the $2$-graph $\Lambda(T,q,t,w)$ associated to the parameters $(T,q,t,w)$ consist of functions $v : T \to \Z / q \Z$ such that $\sum_{i \in T} v(i)w(i)  = t \pmod q$. The parameters $T = \{ 0, e_1, e_2, e_1 + e_2 \}$, $q = 2$, $t = 0$ and $w \equiv 0$ are the only ones which yield the eight vertices shown in \eqref{eq:8tiles}.  However $\Lambda (T,2,0,0)$ has sixteen vertices, hence the $2$-graph $\Lambda_{BD}$ of Example \ref{ldr} (2) is \textbf{not} in the family of $2$-graphs described in \cite{prw}.
\end{rem}

\noindent It turns out that all 2-graphs described in \cite{prw} can be replicated using the basic data given in this paper.



\begin{prop} \label{bdreln}
Let $T$ be a nondegenerate tile, $q \ge 2$ an integer, $t \in \Z / q \Z$ ,and $w: T \to \Z / q \Z$ a function
with $w( c_1 e_1), w( e_2 c_2 )$ invertible in $\Z / q \Z$. Let $\Lambda(T,q,t,w)$ be the associated $2$-graph (see \cite{prw}). For each $p \in ( \Z / q \Z )^P$ there is a bijection $f_p : \Z / q \Z \to \Z / q \Z$ such that
$\Lambda (T,q,t,w) \cong \Lambda_{BD}$ where  $BD = ( T, \Z / q \Z , \{ f_p : p \in ( \Z / q \Z )^P \})$.
\end{prop}

\begin{proof}
For $p \in ( \Z / q \Z )^P$ and $a \in \Z / q \Z$ we may define
$f_p : \Z / q \Z \to \Z / q \Z$ by
\begin{equation} \label{eq:fpformula0}
f_p ( a ) = \Big( t - \sum_{i \in P} w(i) p(i) - w(c_2 e_2) a \Big) w(c_1 e_1)^{-1}
\end{equation}

\noindent since $w(c_1e_1)$ is invertible in $\Z / q \Z$. Observe that
\begin{equation} \label{eq:fpformula}
\sum_{i \in P} w(i) p(i) + w(c_1 e_1) f_p (a) + w( c_2 e_2 ) a = t \text{ for all $a \in \Z / q \Z$.}
\end{equation}

\noindent To see that
$f_p : \Z / q \Z \to \Z / q \Z$ is a bijection we solve \eqref{eq:fpformula} for $a$, and then define
\[
f_p^{-1} (a) = \Big( t - \sum_{i \in P} w(i) p(i) - w(c_1 e_1) f_p (a) \Big) w(c_2 e_2)^{-1} .
\]

\noindent Which makes sense since $w(c_2e_2)$ is invertible in $\Z / q \Z$.

Set $BD = ( T , \Z / q \Z , \{ f_p : p \in ( \Z / q \Z )^P \} )$ and let $v \in \Lambda(T,q,t,w)^0$. We define $p(v) \in ( \Z / q \Z )^P$ by
$p(v) = v |_P$. Define the map $\phi : \Lambda ( T,q,t,w)^0 \to \Lambda_{BD}^0$
by $\phi (v) = F_{p(v),v(c_2 e_2)}$ where $p(v) : \Z / q \Z \to \Z / q \Z$ is the bijection defined in \eqref{eq:fpformula0}. For $F_{p,a} \in \Lambda_{BD}^0$ define $v_{p,a} : T \to \Z / q \Z$ by $v_{p,a} |_P = p$, $v_{p,a} (c_2e_2)=a$ and $v_{p,a} (c_1e_1)=f_p(a)$. By \eqref{eq:fpformula} we see that $v_{p,a} \in \Lambda (T,q,t,w)^0$. Hence $F_{p,a} \mapsto v_{p,a}$ defines a map $\psi : \Lambda_{BD}^0 \to \Lambda(T,q,t,w)^0$ such that $\phi \circ \psi ( F_{p,a} ) = F_{p,a}$ for all $p \in ( \Z / q \Z )^P$ and $a \in \Z / q \Z$. It follows that $\Lambda (T,q,t,w)^0 = \Lambda_{BD}^0$ as functions from $T$ to $\Z / q \Z$. Similarly, we have $\Lambda (T,q,t,w)^n = \Lambda_{BD}^n$ as maps from $T(n)$ to $\Z / q \Z$ for all $n \in \mathbb{N}^2$.
Hence $\Lambda (T,q,t,w)^* = \Lambda_{BD}^*$ and since the structure maps and the degree map are defined in the same manner it follows that the $2$-graphs $\Lambda (T,q,t,w)$ and $\Lambda_{BD}$ are isomorphic, as required.
\end{proof}


\begin{rem} \label{rem:T01} If $T=\{0\}$ is a degenerate tile, then one checks that $\Lambda (T,q,t,w) \cong \Lambda_{BD}$ where
$BD=(T,\Z / q \Z , \{ f_{t w(0)^{-1}} \} )$.
\end{rem}

        \noindent We now turn to the proof of Theorem \ref{2graph}. Our goal is for $\Lambda_{BD}^*$ to be the morphisms in a category so that from a set of basic data we can generate a $2$-graph. To this end we need to say how compose elements $\mu, \nu \in \Lambda_{BD}^*$. Our method is motivated by \cite{prw}.

        \begin{prop}\label{comp}

        Suppose that $\mu, \nu \in \Lambda_{BD}^*$ satisfy $s(\mu) = r(\nu)$. Then there exists a unique path $\lambda \in \Lambda_{BD}^{d(\mu) + d(\nu)}$ such that
        \begin{align}
                \lambda \big( 0, d(\mu) \big) = \mu ~ \text{and} ~ \lambda \big( d(\mu), d(\mu) + d(\nu) \big) = \nu. \label{compeq}
        \end{align}
        \end{prop}

        Notice that $\lambda$ is already explicitly defined on $T \big( d(\mu) \big) \cup \big( T \big( d(\nu) \big) + d(\mu) \big)$. Since $\mu \nu \in \Lambda_{BD}^*$ it remains to show that there is a unique function $\lambda'$ which extends $\lambda$ to $T \big( d(\mu) + d(\nu) \big)$ in such a way that $\lambda'\addon_{T + k} \in \Lambda_{BD}^0$ for every $k \in T \big( d(\mu) + d(\nu) \big) \setminus \Big( T \big( d(\mu) \big) \cup \big( T \big( d(\nu) \big) + d(\mu) \big) \Big)$

        In \cite[Lemma~3.3]{prw} we are given an algorithm for extending $\lambda$ from $T \big( d(\mu) \big) \cup \big( T \big( d(\nu) \big) + d(\mu) \big)$ to $T \big( d(\mu) + d(\nu) \big)$. This involves filling in the missing regions one square at a time so that at each stage there is only one unique choice for the additional value we wish to add. The next two lemmas are an analogue of this.

    \begin{lemma}\label{add}

    Let $T$ be a nondegenerate tile. Fix $n \in \N^2$ and $S$ a subset of $\N^2$ containing $T + n + e_1$ and $T + n - e_2$. Suppose that $\lambda : S \rightarrow A$ is a function such that $\lambda\addon_{T + k}$ belongs to $\Lambda_{BD}^0$ for every $k \in \N^2$ such that $T + k \subset S$. Then there is a unique function
    \[
            \lambda' : S' := S \cup \{ n + (c_1 + 1)e_1 - e_2 \} \rightarrow A
    \]
such that $\lambda'|_S = \lambda$ and $\lambda'\addon_{T + k}$ belongs to $\Lambda_{BD}^0$ for every $k \in \N^2$ such that $T + k \subset S'$.

        \end{lemma}
        \begin{proof}

        The proof here closely follows the proof of \cite[Lemma~3.3]{prw}. Fix $n \in \N^2$ and $S \subset \N^2$ as in the statement of the Lemma. If $n + (c_1 + 1)e_1 - e_2 \in S$ then there is nothing to do. So suppose $n + (c_1 + 1)e_1 - e_2 \notin S$. Let $i \in T \setminus \{c_1 e_1 \}$. Then either $i_2 = 0$ and $i_1 < c_1$ or $i_2 > 0$. If $i_2=0$ and $i_1 < c_1$ we have $i + e_1 \in T$ and hence $n + e_1 - e_2 + i \in T + n - e_2$. Alternatively if $i_2 > 0$, we have $i - e_2 \in T$ and it follows that $n + e_1 - e_2 + i \in T + n + e_1$. Hence $(T \setminus \{ c_1 e_1 \}) + n + e_1 - e_2 \subseteq S$ so the corner $n + (c_1 + 1)e_1 - e_2$ is the only element of $T + n + e_1 - e_2$ which is not in $S$. To say what that point is we now deviate from \cite[Proof of Lemma~3.3]{prw}.

        Let $T' = (T \setminus \{ c_1 e_1 \}) + n + e_1 - e_2 \subseteq S$, $p = \lambda|_{T' \setminus \{ c_2 e_2 \}}$, and $t := \big( n + e_1 + (c_2 - 1)e_2 \big) \in T'$. Then $t$ is the top left hand corner of $T'$. We now define $\lambda'|_S = \lambda$ and
        \begin{align}
                \lambda' \big( n + (c_1 + 1)e_1 - e_2 \big) &:= F_{p, \lambda(t)} (c_1 e_1) = f_p \big( \lambda(t) \big). \label{rule}
        \end{align}

        To see that this is the required function fix $k \in \N^2$ such that $T \setminus \{c_1 e_1\} \subset S'$. If $T + k \subset S$, we have $\lambda'\addon_{T + k} = \lambda\addon_{T + k} \in \Lambda_{BD}^0$. On the other hand, if $T + k \not\subset S$ we must have $k = n + e_1 - e_2$. Then $\lambda'|_{T \setminus \{ c_1 e_1 \} + k} = \lambda|_{T \setminus \{ c_1 e_1 \} + k}$ and (\ref{rule}) gives the unique value $\lambda'(n + e_1 - e_2 + c_1 e_1)$ such that $\lambda'\addon_{T + k} \in \Lambda_{BD}^0$ as required.
        \end{proof}

        \begin{lemma}\label{add2}

        Let $T$ be a nondegenerate tile. Fix $n \in \N^2$ and $S$ a subset of $\N^2$ containing $T + n + e_1$ and $T + n - e_2$. Suppose that $\lambda : S \rightarrow A$ is a function such that $\lambda\addon_{T + k}$ belongs to $\Lambda_{BD}^0$ for every $k \in \N^2$ such that $T + k \subset S$. Then there is a unique function
    \[
            \lambda' : S' := S \cup \{ n + c_2 e_2 \} \rightarrow A
    \]
such that $\lambda'|_S = \lambda$ and $\lambda'\addon_{T + k}$ belongs to $\Lambda_{BD}^0$ for every $k \in \N^2$ such that $T + k \subset S'$.

        \end{lemma}

        \begin{remark*} The statement of Lemma~\ref{add2} is very similar to that of Lemma~\ref{add} and so we omit it. The key point of difference is that for Lemma~\ref{add2} we are filling in the upper left hand corner $\lambda' (n + c_2 e_2)$ while in Lemma~\ref{add} we filled in the lower right hand corner $\lambda' (n + (c_1 + 1)e_1 - e_2)$. A mirror argument of the proof of Lemma~\ref{add} using the fact that $f_p$ is a bijection for each $p \in A^P$ will prove Lemma \ref{add2}.
        \end{remark*}

        \begin{proof}[Proof of Proposition~\ref{comp}]

        To obtain $\lambda \in \Lambda_{BD}^{d ( \mu ) + d ( \nu )}$ we start with the region $T \big( d(\mu) \big) \cup \big( T \big( d(\nu) \big) + d(\mu) \big)$. We must fill in the two remaining rectangular regions given by
        \begin{align*}
                BR &:= \{ j \in \N^2 : c_1 + d(\mu)_1 < j_1 \leq c_1 + d(\mu)_1 + d(\nu)_1, j_2 < d(\mu)_2 \}, ~ \text{and} \\
                UL &:= \{ j \in \N^2 : j_1 < d(\mu)_1, c_2 + d(\mu)_2 < j_2 \leq c_2 + d(\mu)_2 + d(\nu)_2 \}.
        \end{align*}
We consider each region separately. First consider the bottom right rectangle $BR$. We order $BR$ lexicographically by starting at the top of the $(c_1 + d(\mu)_1 + 1)$th column and going down. We then repeat with the $(c_1 + d(\mu)_1 + 2)$th column and so on. Now we proceed by induction. For our base case consider $d(\mu) + (c_1 + 1)e_1 - e_2 \in BR$. We apply Lemma~\ref{add} with $n = d(\mu)$ to give $\lambda (d(\mu) + (c_1 + 1)e_1 - e_2) = f_p \big( \lambda \big( d(\mu) + c_1 e_1 - e_2 \big) \big)$ for $p = \lambda|_{P + d(\mu) - e_2}$.

        Fix $k \in BR$. Suppose $\lambda'\addon_{T + j} \in \Lambda_{BD}^0$ for every $j$ such that $d(\mu) + (c_1 + 1)e_1 - e_2 \leq j \leq k$ in the lexicographic order described above. For the inductive step there are two cases to consider. Firstly, if $k_2 \neq 0$, let $n = k - (c_1 + 1)e_1$. Then $\lambda'(n + e_1)$ and $\lambda'(n - e_2)$ are given by the inductive hypothesis. Then Lemma~\ref{add} with this $n$ gives the unique value for $\lambda'(k - e_2)$. On the other hand, if $k_2 = 0$ then the next value we need to determine is $\lambda \big( k + e_1 + (d(\mu)_2 - 1)e_2 \big)$, that is, the top value of the next column in $BR$. Then $\lambda \big( T + k + d(\mu)_2 e_2 \big) = \nu \big( T + (k_1 - d(\mu)_1)e_1 \big)$ and $T \big( k - e_1 + (d(\mu)_2 - 1)e_2 \big)$ is given by the inductive hypothesis. Thus Lemma~\ref{add} with $n = k - c_1 e_1 + d(\mu)_2 e_2$ gives the unique value $\lambda \big( k + e_1 + (d(\mu)_2 - 1)e_2 \big)$.

        It remains to show how to fill in the upper left rectangle $UL$. We again order elements lexicographically. However this time we start at the right hand end of the $\big( c_2 + d(\mu)_2 + 1 \big)$th row and count left. At the end of that row we then continue from the right hand end of the $\big( c_2 + d(\mu)_2 + 2 \big)$th row, and so on. For a base case for an induction argument we consider $d(\mu) - e_1 + e_2$. A mirror argument of the argument in the preceding two paragraphs, replacing the use of Lemma~\ref{add} with Lemma~\ref{add2} will complete the proof.
        \end{proof}

    \noindent We now know how to compose elements $\mu, \nu \in \Lambda_{BD}^*$ and can hence prove Theorem~\ref{2graph}.

    \begin{proof}[Proof of Theorem~\ref{2graph}]

        Proposition~\ref{comp} ensures that composition in $(\Lambda_{BD}^0, \Lambda_{BD}^*, r, s)$ is well defined. We now show that composition in $\Lambda_{BD}$ is associative. Fix $\mu, \nu, \rho \in \Lambda_{BD}^*$ such that $s(\mu) = r(\nu)$ and $s(\nu) = r(\rho)$. Then Proposition~\ref{comp} implies that $\mu \nu$ is the unique path in $\Lambda_{BD}^{d(\mu) + d(\nu)}$ such that $(\mu \nu) \big( 0, d(\mu) \big) = \mu$ and $(\mu \nu) \big( d(\mu), d(\mu) + d(\nu) \big) = \nu$. Proposition~\ref{comp} also implies that $\nu \rho$ is the unique path in $\Lambda_{BD}^{d(\nu) + d(\rho)}$ such that $(\nu \rho) \big( 0, d(\nu) \big) = \nu$ and $(\nu \rho) \big( d(\nu, d(\nu) + d(\rho) \big) = \rho$. Then
        \begin{align*}
                (\mu \nu) \rho &= \Big( (\mu \nu) \big( 0, d(\mu) \big) \Big) \Big( (\mu \nu) \big( d(\mu, d(\mu) + d(\nu) \big) \Big) \rho
                        = \mu \nu \rho \\
                        &= \mu \Big( (\nu \rho) \big( 0, d(\nu) \big) \Big) \Big( (\nu \rho) \big( d(\nu), d(\nu) + d(\rho) \big) \Big)
                        = \mu (\nu \rho).
        \end{align*}

        \noindent Hence $\Lambda_{BD}$ is a category, which is countable as Lemma~\ref{add} and Lemma~\ref{add2} imply that $|\Lambda_{BD}^{e_1} \cup \Lambda_{BD}^{e_2} | < \infty$. The map $d : \Lambda_{BD} \rightarrow \N^2$ satisfies $d(\mu \nu) = d(\mu) + d(\nu)$ and so it is a functor. To see that $\Lambda_{BD}$ is a 2-graph it remains to show that $d$ satisfies the factorisation property. Fix $m, n \in \N^2$ and $\lambda \in \Lambda_{BD}^{m + n}$. Then (\ref{compeq}) ensures $\mu := \lambda(0, m)$ and $\nu := \lambda(m, m + n)$ are the two unique paths with $\lambda = \mu \nu$.
        \end{proof}

        \begin{remark} \label{rem:T02}
        If $|A| = 1$ then $\Lambda_{BD}$ is the $2$-graph consisting of just one vertex with one blue loop and one red loop. If $T=\{ 0 \}$
        then every vertex in $\Lambda_{BD}$ connects to every other vertex in $\Lambda_{BD}$ with exactly one blue and one red edge.
        \end{remark}

        \begin{example}

        The $2$-graph associated to the basic data given in Example~\ref{ldr} (2) has the skeleton given below.
        
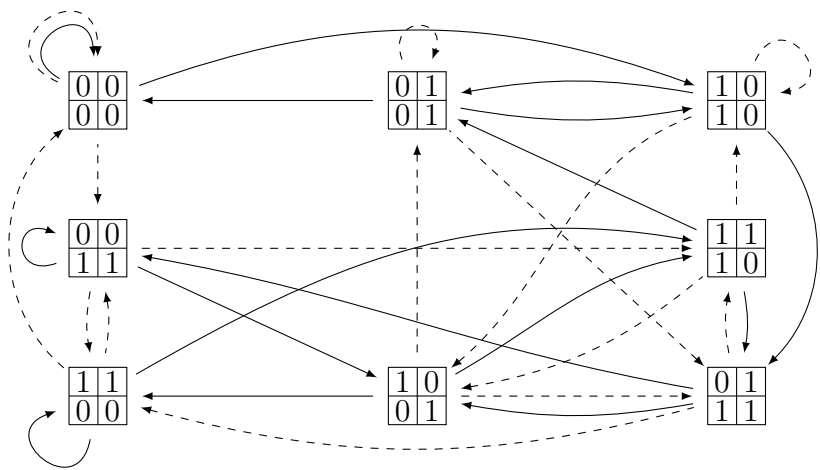
\[
        \begin{tikzpicture}[scale=0.28]
                \node[circle,inner sep=0.5pt] (0000) at (-15, 7) {\begin{tikzpicture}[scale=0.38]
                                                                                                                \draw (0, 0) -- (2, 0) -- (2, 2) -- (0, 2) -- (0, 0);
                                                                                                                \draw (1, 0) -- (1, 2);
                                                                                                                \draw (0, 1) -- (2, 1);
                                                                                                                \node[circle,inner sep=0.5pt] () at (0.5, 0.5) {$0$};
                                                                                                                \node[circle,inner sep=0.5pt] () at (1.5, 0.5) {$0$};
                                                                                                                \node[circle,inner sep=0.5pt] () at (0.5, 1.5) {$0$};
                                                                                                                \node[circle,inner sep=0.5pt] () at (1.5, 1.5) {$0$};
                                                                                                    \end{tikzpicture}}
                        edge[-latex, loop, out=150, in=95, min distance=25, looseness=4] (0000)
                    edge[-latex, loop, out=155, in=90, min distance=35, dashed, looseness=5] (0000)
                        ;
            \node[circle,inner sep=0.5pt] (0101) at (0, 7) {\begin{tikzpicture}[scale=0.38]
                                                                                                                \draw (0, 0) -- (2, 0) -- (2, 2) -- (0, 2) -- (0, 0);
                                                                                                                \draw (1, 0) -- (1, 2);
                                                                                                                \draw (0, 1) -- (2, 1);
                                                                                                                \node[circle,inner sep=0.5pt] () at (0.5, 0.5) {$0$};
                                                                                                                \node[circle,inner sep=0.5pt] () at (1.5, 0.5) {$1$};
                                                                                                                \node[circle,inner sep=0.5pt] () at (0.5, 1.5) {$0$};
                                                                                                                \node[circle,inner sep=0.5pt] () at (1.5, 1.5) {$1$};
                                                                                                    \end{tikzpicture}}
                        edge[-latex] (0000)
                        edge[-latex, dashed, loop, out=110, in=70, looseness=4] (0101)
                        ;
            \node[circle,inner sep=0.5pt] (1010) at (15, 7) {\begin{tikzpicture}[scale=0.38]
                                                                                                                \draw (0, 0) -- (2, 0) -- (2, 2) -- (0, 2) -- (0, 0);
                                                                                                                \draw (1, 0) -- (1, 2);
                                                                                                                \draw (0, 1) -- (2, 1);
                                                                                                                \node[circle,inner sep=0.5pt] () at (0.5, 0.5) {$1$};
                                                                                                                \node[circle,inner sep=0.5pt] () at (1.5, 0.5) {$0$};
                                                                                                                \node[circle,inner sep=0.5pt] () at (0.5, 1.5) {$1$};
                                                                                                                \node[circle,inner sep=0.5pt] () at (1.5, 1.5) {$0$};
                                                                                                    \end{tikzpicture}}
                        edge[-latex, out=170, in=10] (0101)
                        edge[latex-, out=-170, in=-10] (0101)
                        edge[-latex, dashed, loop, out=65, in=10, min distance=25, looseness=4] (1010)
                        edge[latex-, out=160, in=20] (0000)
                        ;
            \node[circle,inner sep=0.5pt] (1011) at (15, 0) {\begin{tikzpicture}[scale=0.38]
                                                                                                                \draw (0, 0) -- (2, 0) -- (2, 2) -- (0, 2) -- (0, 0);
                                                                                                                \draw (1, 0) -- (1, 2);
                                                                                                                \draw (0, 1) -- (2, 1);
                                                                                                                \node[circle,inner sep=0.5pt] () at (0.5, 0.5) {$1$};
                                                                                                                \node[circle,inner sep=0.5pt] () at (1.5, 0.5) {$0$};
                                                                                                                \node[circle,inner sep=0.5pt] () at (0.5, 1.5) {$1$};
                                                                                                                \node[circle,inner sep=0.5pt] () at (1.5, 1.5) {$1$};
                                                                                                    \end{tikzpicture}}
                        edge[-latex, dashed] (1010)
                        edge[-latex] (0101)
                        ;
            \node[circle,inner sep=0.5pt] (1101) at (15, -7) {\begin{tikzpicture}[scale=0.38]
                                                                                                                \draw (0, 0) -- (2, 0) -- (2, 2) -- (0, 2) -- (0, 0);
                                                                                                                \draw (1, 0) -- (1, 2);
                                                                                                                \draw (0, 1) -- (2, 1);
                                                                                                                \node[circle,inner sep=0.5pt] () at (0.5, 0.5) {$1$};
                                                                                                                \node[circle,inner sep=0.5pt] () at (1.5, 0.5) {$1$};
                                                                                                                \node[circle,inner sep=0.5pt] () at (0.5, 1.5) {$0$};
                                                                                                                \node[circle,inner sep=0.5pt] () at (1.5, 1.5) {$1$};
                                                                                                    \end{tikzpicture}}
                    edge[-latex, dashed, out=100, in=-100] (1011)
                        edge[latex-, out=80, in=-80] (1011)
                        edge[latex-, out=45, in=-45] (1010)
                        edge[latex-, dashed] (0101)
                ;
            \node[circle,inner sep=0.5pt] (0110) at (0, -7) {\begin{tikzpicture}[scale=0.38]
                                                                                                                \draw (0, 0) -- (2, 0) -- (2, 2) -- (0, 2) -- (0, 0);
                                                                                                                \draw (1, 0) -- (1, 2);
                                                                                                                \draw (0, 1) -- (2, 1);
                                                                                                                \node[circle,inner sep=0.5pt] () at (0.5, 0.5) {$0$};
                                                                                                                \node[circle,inner sep=0.5pt] () at (1.5, 0.5) {$1$};
                                                                                                                \node[circle,inner sep=0.5pt] () at (0.5, 1.5) {$1$};
                                                                                                                \node[circle,inner sep=0.5pt] () at (1.5, 1.5) {$0$};
                                                                                                    \end{tikzpicture}}
                    edge[-latex, dashed] (1101)
                        edge[latex-, out=-10, in=-170] (1101)
                        edge[-latex, dashed] (0101)
                        edge[latex-, dashed, out=40, in=-160] (1010)
                        edge[latex-, dashed, out=10, in=-140] (1011)
                        edge[-latex, out=30, in=-170] (1011)
                    ;
            \node[circle,inner sep=0.5pt] (0011) at (-15, -7) {\begin{tikzpicture}[scale=0.38]
                                                                                                                \draw (0, 0) -- (2, 0) -- (2, 2) -- (0, 2) -- (0, 0);
                                                                                                                \draw (1, 0) -- (1, 2);
                                                                                                                \draw (0, 1) -- (2, 1);
                                                                                                                \node[circle,inner sep=0.5pt] () at (0.5, 0.5) {$0$};
                                                                                                                \node[circle,inner sep=0.5pt] () at (1.5, 0.5) {$0$};
                                                                                                                \node[circle,inner sep=0.5pt] () at (0.5, 1.5) {$1$};
                                                                                                                \node[circle,inner sep=0.5pt] () at (1.5, 1.5) {$1$};
                                                                                                    \end{tikzpicture}}
                        edge[latex-] (0110)
                        edge[latex-, dashed, out=-15, in=-165] (1101)
                        edge[-latex, out=30, in=170] (1011)
                        edge[-latex, loop, out=-100, in=-160, looseness=4] (0011)
                        edge[-latex, dashed, out=140, in=-140] (0000)
                        ;
                \node[circle,inner sep=0.5pt] (1100) at (-15, 0) {\begin{tikzpicture}[scale=0.38]
                                                                                                                \draw (0, 0) -- (2, 0) -- (2, 2) -- (0, 2) -- (0, 0);
                                                                                                                \draw (1, 0) -- (1, 2);
                                                                                                                \draw (0, 1) -- (2, 1);
                                                                                                                \node[circle,inner sep=0.5pt] () at (0.5, 0.5) {$1$};
                                                                                                                \node[circle,inner sep=0.5pt] () at (1.5, 0.5) {$1$};
                                                                                                                \node[circle,inner sep=0.5pt] () at (0.5, 1.5) {$0$};
                                                                                                                \node[circle,inner sep=0.5pt] () at (1.5, 1.5) {$0$};
                                                                                                    \end{tikzpicture}}
                        edge[latex-, dashed] (0000)
                        edge[-latex, dashed] (1011)
                        edge[-latex] (0110)
                        edge[latex-, dashed, out=-80, in=80] (0011)
                        edge[-latex, dashed, out=-100, in=100] (0011)
                        edge[-latex, loop, out=-160, in=160, looseness=4] (1100)
                        edge[latex-, out=-10, in=170] (1101)
                        ;
    \end{tikzpicture}
        \]



       \end{example}

       \noindent Next we prove some general facts about the skeleton of a basic data 2-graphs.

        \begin{prop}\label{sk}

        Let $\Lambda_{BD}$ be the $2$-graph associated to the basic data $BD$, then

        \begin{enumerate}

        \item $|\Lambda_{BD}^0| = |A|^{|P| + 1}$;

        \item for any $u, v \in \Lambda_{BD}^0$ and $i \in \{ 1, 2 \}$ we have $v \Lambda_{BD}^{e_i} u \neq \emptyset$ if and only if
        \begin{align}
                v(m) = u(m - e_i) ~ \text{for every} ~ m \in T \cap (T + e_i), \label{connect}
        \end{align}
in which case $|v \Lambda_{BD}^{e_i} u| = 1$; and

        \item $|v \Lambda_{BD}^{e_1}| = |A|^{c_2} = |\Lambda_{BD}^{e_1} v|$ and $|v \Lambda_{BD}^{e_2}| = |A|^{c_1} = |\Lambda_{BD}^{e_2} v|$ for every $v \in \Lambda_{BD}^0$.

            \item $| v \Lambda_{BD}^{e_1+e_2} u | =1$ for all $u,v \in \Lambda_{BD}^0$.

        \end{enumerate}

        \end{prop}

        \begin{rem}\label{rfnss}

        Proposition~\ref{sk}~(3) shows that $\Lambda_{BD}$ is row-finite and has no sinks or sources.

        \end{rem}

        \begin{proof}[Proof of Proposition~\ref{sk}]


        \textbf{(1)}: Recall that $\Lambda_{BD}^0 = \{ F_{p, a} : p \in A^P, a \in A \}$. For each $p \in A^P$ and each $a \in A$ the functions $F_{p, a}$ are distinct; hence we have that $| \Lambda_{BD}^0 | = |A^P||A| = |A|^{|P| + 1}$.

        \textbf{(2):} Let $u, v \in \Lambda_{BD}^0$ and $i \in \{ 1, 2 \}$. Suppose $v \Lambda^{e_i} u \neq \emptyset$ and fix $m \in T \cap (T + e_i)$. Then there exists $\lambda \in v \Lambda_{BD}^{e_i} u$. To show \eqref{connect} holds we calculate
        \[
                v(m) = \lambda\addon_T (m) = \lambda|_T (m) = \lambda|_{T + e_i} (m - e_i) = u (m - e_i).
        \]
For the reverse implication suppose $u, v$ and $i$ satisfy (\ref{connect}). Define $\lambda : T(e_i) \rightarrow A$ by
        \begin{align*}
                \lambda (m) &= \begin{cases}
                                            v(m) &\text{for} ~ m \in T \\
                                            u(m) &\text{for} ~ m \in (T + e_i) \setminus T.
                                            \end{cases}
        \end{align*}
Then $r(\lambda) = \lambda\addon_T = v$ and \eqref{connect} ensures that for every $m \in T \cap (T + e_i)$ we have $\lambda (m) = v(m) = u (m - e_i)$, which forces $s(\lambda) = \lambda\addon_{T + e_i} = u$. Hence $\lambda \in v \Lambda_{BD}^{e_i} u$. To see that $|v \Lambda_{BD}^{e_i} u| = 1$ suppose $\lambda' \in v \Lambda_{BD}^{e_i} u$. Then $s(\lambda') = u$, $r(\lambda') = v$ and \eqref{connect} forces $\lambda' = \lambda$.

        \textbf{(3):} Fix $v \in \Lambda_{BD}^0$ and let $S = \{ u \in \Lambda_{BD}^0 : v \Lambda^{e_1}_{BD} u \neq \emptyset \}$. Then Proposition~\ref{sk}~(2) implies $|v\Lambda^{e_1}_{BD}| = |v\Lambda^{e_1}_{BD} S| = |S|$. We claim $|S| = |A|^{c_2}$. For this we construct $u \in \Lambda^0_{BD}$ such that $u \in S$. Let $u|_{T \cap (T - e_1)} = v|_{T \cap (T + e_1)}$. Then there remains $|T| - |T \cap (T - e_1)| = c_2 + 1$ values left to be determined, one of which must be the corner $u(c_1 e_1)$. We may arbitrarily define values for $u|_{(T \cap (T - e_1)) \setminus \{ c_1 e_1 \}}$, which may be done in $|A|^{c_2}$ ways. Then $u(c_1 e_1) = F_{u|_P, u(c_2 e_2)} (c_1 e_1) = f_{u|_P} (c_2 e_2)$. Hence $|v \Lambda_{BD}^{e_1}| = |U| = |A|^{c_2}$.

        To see that $|\Lambda_{BD}^{e_1} v| = |A|^{c_2}$ let $R = \{ u \in \Lambda_{BD}^0 : u \Lambda_{BD}^{e_1} v \neq \emptyset \}$. Then a similar argument to that given above shows that $|\Lambda_{BD}^{e_1} v| = |R| = |A|^{c_2}$. The mirror results for edges of degree $e_2$ follow by symmetry.

        \textbf{(4):} Follows by applying Lemma \ref{add} and Lemma \ref{add2} to $S = T \cup ( T + e_1+e_2 )$.
        \end{proof}

%
%

\begin{defn} A $k$-graph $\Lambda$ is \emph{strongly connected} if for every $u, v \in \Lambda^0$ there exists $\lambda \in \Lambda$ such that $r(\lambda) = v$ and $s(\lambda) = u$.
\end{defn}

\noindent It turns out that $\Lambda_{BD}$ is always strongly connected.

        \begin{prop}\label{sc}

        Let $\Lambda_{BD}$ be the $2$-graph associated to the basic data $BD$. Suppose $k \in \N$ satisfies $(k - 1)(e_1 + e_2) \in T$ and $k(e_1 + e_2) \notin T$. Then for every $u, v \in \Lambda_{BD}^0$ there exists $\lambda \in \Lambda_{BD}^{k(e_1 + e_2)}$ such that $r(\lambda) = v$ and $s(\lambda) = u$. Hence $\Lambda_{BD}$ is strongly connected.

        \end{prop}
        \begin{proof}

        We will use the following notation: For $n \in \N$ let $T_n := T + n(e_1 + e_2)$, $P_n := P + n(e_1 + e_2)$ and $\mu^n \in \Lambda_{BD}^{n(e_1 + e_2)}$. Fix $u, v \in \Lambda^0$ and $\mu^m \in v \Lambda^{m(e_1 + e_2)}$. Then
        \begin{align}
                s(\mu^m)(i) = u \big( i - (k - m)(e_1 + e_2) \big) ~ \text{for every} ~ i \in T \cap T_{k - m}. \label{A}
        \end{align}
We show that $\mu^m$ exists for $m \leq k$ so that $\mu^k$ will be the required path. We do this by an induction argument. For a base case let $\mu^0 = v$. Now suppose for $0 \leq m < k$ that $s(\mu^m)$ satisfies (\ref{A}). To show that $s(\mu^{m+1})$ satisfies (\ref{A}) it suffices to show there exists $\lambda \in \Lambda_{BD}^{e_1 + e_2}$ such that $r(\lambda) = s(\mu^m)$ and
        \begin{align}
                s(\lambda)(i) = u \big( i - \big( k - (m + 1) \big) (e_1 + e_2) \big) ~ \text{for} ~ i \in T \cap T_{k - (m + 1)} \label{B}
        \end{align}
so that Proposition~\ref{comp} implies $\mu^m \lambda = \mu^{m + 1}$.

        Let $S = T_m \cup T_{m + 1}$. We will define a function $\lambda' : S \rightarrow A$. Let $\lambda'\addon_{T_m} = s(\mu^m)$ and for $i \in T_k \cap T_{k - (m + 1)}$ let $\lambda' (i) = u \big( i - k(e_1 + e_2) \big)$. This fixes those values of $\lambda'$ whose domain coincides with the domain of $s(\mu^m)$ or the domain of $u$. We may now define the values of $\lambda'|_{P_{m+1}}$ and $j = \lambda' \big( m(e_1 + e_2) + c_2 e_2 \big)$ arbitrarily. Then $\lambda' \big( m(e_1 + e_2) + c_1 e_1 \big) = f_{\lambda'|_{P_{m + 1}}} (j)$. This completes the definition of $\lambda'$. Now Lemma~\ref{add} and Lemma~\ref{add2} with these $S$ and $\lambda'$ give $\lambda \in \Lambda_{BD}^{e_1 + e_2}$ such that (\ref{B}) holds and $r(\lambda) = \lambda'\addon_{T_m} = s(\mu^m)$. The final remark now follows from the definition.
        \end{proof}

        \section{Aperiodicity}\label{aps}

        The results of \cite{kp,ls, rsy1, shotwell} show that a condition, called aperiodicity, plays an important role
        in simplicity results for the $C^*$-algebras associated to $k$-graphs. In this section, we prove Theorem~\ref{APT}
        which gives conditions on the
        basic data $BD$ which ensure that the associated $2$-graph $\Lambda_{BD}$ is aperiodic. Theorem~\ref{APT} generalises
       \cite[Theorem 5.2]{prw} as it allows us to prove aperiodicity for a larger class of $2$-graphs.

        There have been many formulations of aperiodicity, see \cite{fmy, kp, ls, rsy1, robsi1, shotwell}. The new formulation we propose is inspired by condition (4) given by Robertson and Sims in~\cite[Lemma 3.3]{robsi1}.

 \begin{defn}
 Let $\Lambda$ be a row-finite $k$-graph with no sources then $\Lambda$ is \emph{aperiodic}
 if for every $v \in \Lambda^0$ and for every distinct $m, n \in \N^k$ with $m \wedge n = 0$ there exists a path $\lambda \in \Lambda$ such that $r(\lambda) = v$, $d(\lambda) \geq m \vee n$ and
        \begin{align}
                \lambda \big( m, m + d(\lambda) - (m \vee n) \big) \neq \lambda \big( n, n + d(\lambda) - (m \vee n) \big). \label{APC}
        \end{align}
\end{defn}

\noindent The following lemma shows that our definition of aperiodicity is equivalent to condition (4) of \cite[Lemma 3.3]{robsi1}. Briefly, condition (4) of \cite[Lemma 3.3]{robsi1}, which we shall refer to as the \emph{Robertson-Sims condition} amounts to demanding that \eqref{APC} holds for all $m,n \in \mathbb{N}^2$.

\begin{lemma} \label{lem:reducecase}
Let $\Lambda$ be a row-finite $2$-graph with no sources. Then $\Lambda$ is aperiodic if and only if it satisfies the Robertson-Sims condition.
\end{lemma}

\begin{proof}
Suppose that $\Lambda$ is aperiodic. Fix $v \in \Lambda^0$ and suppose that $m,n \in \mathbb{N}^2$ are such that $m \wedge n \neq 0$. Put $m' = m - (m \wedge n)$ and $n' = n - ( m \wedge n)$. Then $m' \wedge n'=0$ and $m' \vee n' = ( m \vee n ) - ( m \wedge n)$. Let $w$ be such that $v \Lambda^{m \wedge n} w \neq \emptyset$. Then by hypothesis there exists $\lambda' \in w \Lambda$ with $d ( \lambda' ) \ge m' \vee n'$ and
\[
        \lambda' \big(m'  , m' + d(\lambda') - (m' \vee n' ) \big) \neq \lambda' \big( n' , n' + d(\lambda') - (m' \vee n' ) \big) .
\]
Let $\mu \in v \Lambda^{m \wedge n} w$. Then $\lambda = \mu \lambda'$ is such that
\begin{align*}
\lambda \big( m, m + d(\lambda) - (m \vee n) \big) &= \lambda' \big( m - (m \wedge n), m - (m \wedge n) + d(\lambda) - (m \vee n) \big) \\
        &= \lambda' \big( m', m' + d(\lambda') + (m \wedge n ) - (m \vee n) \big) \\
        &= \lambda' \big( m', m' + d(\lambda') - (m' \vee n' ) \big) \\
                &\neq \lambda' \big( n', n' + d(\lambda') - (m' \vee n') \big) \\
        &= \lambda' \big( n', n' + d(\lambda') + (m \wedge n ) - (m \vee n) \big) \\
        &= \lambda' \big( n - (m \wedge n), n - (m \wedge n) + d(\lambda) - (m \vee n) \big) \\
        & = \lambda \big( n, n + d(\lambda) - (m \vee n) \big) ,
\end{align*}
\noindent and so $\Lambda$ satisfies the Robertson-Sims condition. If $m \wedge n = 0$ then the Robertson-Sims condition follows
immediately from the aperiodicity of $\Lambda$.

If $\Lambda$ satisfies the Robertson-Sims condition then it is clear that it is aperiodic; which completes the proof.
\end{proof}




        \noindent We now give a condition on the basic data $BD$ which, as we shall see in Theorem~\ref{APT}, ensures that $\Lambda_{BD}$ is aperiodic.

        \begin{definition}\label{ac}

         Let $\Lambda_{BD}$ be the 2-graph associated to the basic data $BD$. For $a \in A$ we say that $BD$ has a \emph{blue $a$-breaking cycle} if there exists $k > 1$ and distinct $v_1,\ldots,v_k \in \Lambda_{BD}^0$ such that for every $i \in \{1,\ldots,k\}$ we have $v_i|_{T \cap (T + e_2)} \equiv a$ and either

        \begin{enumerate}

        \item $| v_i \Lambda^{e_1} v_{i + 1} | = 1$ for every $i \in \{1,\ldots,k\}$ where $v_{k+1}=v_1$; or

        \item $| v_i \Lambda^{e_1} v_i | = 1$ for every $i \in \{1,\ldots,k\}$.

        \end{enumerate}

        \noindent We may define what it means for $BD$ to have a \emph{red $a$-breaking cycle} in a similar way. We say that $BD$ has an
        \emph{$a$-breaking cycle} if it has either a blue $a$-breaking cycle or a red $a$-breaking cycle.

        \end{definition}

\begin{examples} \label{rem:b4APT}
\begin{enumerate}

\item The basic data given in Example~\ref{ldr} has a blue $0$-breaking cycle as it satisfies part (2) of Definition \ref{ac} where $v_1 = F_{p[0],0}$
and $v_2 = F_{p[1],0}$. Moreover, it has a blue $1$-breaking cycle as it satisfies part (1) of Definition \ref{ac} where $v_1 = F_{p[0],1}$ and $v_2 = F_{p[1],1}$. One also checks that there are a red $1$-breaking cycle and a red $0$-breaking cycle.

\item The basic data given in Example~\ref{ldr} (2) has a blue $0$-breaking cycle as it satisfies part (2) of Definition \ref{ac} where $v_1 = F_{p[0,0],0}$ and $v_2 = F_{p[1,0],1}$. One also checks that there is a red $0$-breaking cycle, but no red $1$-breaking cycle.

\item Let $A=\{0,1\}$, and $T=\{0,e_1,2e_1,e_2 \}$, so that $P=\{ 0 , e_1 \}$.
For $a,b \in \{ 0,1\}$ let $p[a,b] : \{ 0 , e_1 \} \rightarrow A$ be given by $p[a,b](0)=a$ and $p[a,b](e_1)=b$, so
that $A^P = \{ p[0,0], p[0,1], p[1,0], p[1,1] \}$.
Define a bijection $f_{p[0,0]} : A \to A$ by $f_{p[0,0]}(a)=a$ for all $a \in A$,
and bijections $f_{p[1,0]}, f_{p[1,1]} , f_{p[0,1]} : A \to A$ by $f_{p[1,0]}(a)=f_{p[0,1]}(a)=f_{p[1,1]}(a)=a+1 \pmod 2$ for all $a \in A$ which gives us basic data $BD$.
 The basic data $BD$ has a blue $0$-breaking cycle as it satisfies part (2) of Definition \ref{ac} with $v_1=F_{p[0,0],0}$ and $v_2 = F_{p[1,1],1}$. It also has a blue $1$-breaking cycle as it satisfies part (1) of
Definition \ref{ac} with $v_1 = F_{p[0,0],1}$, $v_2 = F_{p[0,1],0}$ and $v_3 = F_{p[1,0],0}$. One checks that there is a red $1$-breaking cycle
but not a red $0$-breaking cycle.

\end{enumerate}
\end{examples}

\noindent We now give the main result of this paper.

        \begin{theorem}\label{APT}

        Let $\Lambda_{BD}$ be the 2-graph associated to the basic data $BD$. If $BD$ has an $a$-breaking cycle for some $a \in A$ then $\Lambda_{BD}$ is aperiodic.

        \end{theorem}
        \begin{proof}

        Fix $v \in \Lambda_{BD}^0$ and distinct $m, n \in \N^2$. We also fix $\alpha \in v \Lambda_{BD}^{m \vee n}$. If $\alpha|_{T + m} \neq \alpha|_{T + n}$ then $\lambda := \alpha$ will do. So suppose otherwise. Recall from Lemma \ref{lem:reducecase} that it is enough show that \eqref{APC} is satisfied for every distinct $m, n \in \N^2$ such that $m \wedge n = 0$.

				Without loss of generality we may suppose $m_1 > n_1 = 0$ and $n_2 > m_2 = 0$. The proof proceeds as follows: In the case where $BD$ has a blue $a$-breaking cycle we first show there exists $\nu \in \Lambda_{BD}^{m \vee n}$ whose bottom $c_2$ rows are all $a$'s. We then attach an edge $\mu$ to the head of $\nu$ so that $\rho = (\mu \nu)(0, m \vee n)$ satisfies $\rho(m, m) \neq \rho(n, n)$. We then use the fact that $\Lambda_{BD}$ is strongly connected to prove there exist $\beta \in s(\alpha) \Lambda_{BD} r(\rho)$ and then show that $\alpha \beta \rho$ satisfies (\ref{APC}).

        Let $T^- \subset \N^2$ denote the bottom $c_2$ rows of $T$. We claim that for all $l \ge 0$ there exists $\nu' \in \Lambda_{BD}^{l e_1}$ such that $\nu'|_{T^-(l e_1)} \equiv a$. We shall prove our claim by induction: For a base case, $l=0$, fix $p \in A^P$ such that $p \equiv a$. Then $F_{p, f_p^{-1}(a)} \in \Lambda_{BD}^0$ and $F_{p, f_p^{-1}(a)}|_{T^-} \equiv a$. Now suppose that for $k\ge 1$ there is $\phi \in \Lambda_{BD}^{k e_1}$ which satisfies $\phi|_{T^-(k e_1)} \equiv a$. Define $q \in A^P$ by
        \[
                q(j) =
                            \begin{cases}
                                    \phi(j - e_1) & \text{for} ~ j \in P \cap (P + e_1) \\
                                    a & \text{for} ~ j \notin P \cap (P + e_1).
                            \end{cases}
        \]
    Then $F_{q, f_q^{-1}(a)} \in \Lambda_{BD}^0$ and $F_{q, f_q^{-1}(a)}|_{T^-} \equiv a$. Then Proposition~\ref{sk}~(2) implies there exists $\theta \in F_{q, f_q^{-1}(a)} \Lambda_{BD}^{e_1} r(\phi)$. Then $\nu' = \theta \phi \in \Lambda_{BD}^{(k + 1)e_1}$ satisfies $\nu'|_{T^-( \Tiny( k + 1 \Tiny)e_1)} \equiv a$ as required. Hence by the principle of mathematical induction our claim follows.

     Now we may prove the existence of a $\nu \in \Lambda_{BD}^{m \vee n}$ whose bottom $c_2$ rows are all $a$'s. By the claim established in the above paragraph there is $\nu' \in \Lambda_{BD}^{m_1e_1}$ such that $\nu'|_{T^-( m_1 e_1 )} \equiv a$. Since $\Lambda$ has no sources there exists $\nu \in \Lambda_{BD}^{m \vee n}$ such that $\nu(0, m) = \nu'$ and hence $\nu|_{T^-(m)} = \nu'|_{T^-(m)} \equiv a$. In other words, the bottom $c_2$ rows of $\nu$ are all $a$'s.

        By Proposition~\ref{sk}~(2) there exist unique $\mu_1 \in v_1 \Lambda_{BD}^{e_2} r(\nu)$ and $\mu_2 \in v_2 \Lambda^{e_2} r(\nu)$. We claim $(\mu_1 \nu)(m, m) \neq (\mu_2 \nu)(m, m)$. Suppose $\{ v_k \}$ satisfies condition (1) of Definition~\ref{ac}. Then for every $l \leq m$ let $i = 1 + l_1$ (mod k). Then Lemma~\ref{add} implies $(\mu_1 \nu)(l, l) = v_{i} \neq v_{i + 1} = (\mu_2 \nu)(l, l)$. On the other hand, if $\{ v_k \}$ satisfies condition (2) of Definition~\ref{ac} then Lemma~\ref{add} implies $(\mu_1 \nu)(l, l) = v_1 \neq v_2 = (\mu_2 \nu)(l, l)$ for every $l \leq m$. In either case taking $l = m$ proves the claim.

        If $(\mu_1 \nu)(m, m) \neq (\mu_1 \nu)(n, n)$ we let $\rho = (\mu_1 \nu)(0, m \vee n) \in \Lambda_{BD}^{m \vee n}$. Otherwise let $\rho = (\mu_2 \nu)(0, m \vee n) \in \Lambda_{BD}^{m \vee n}$. Then $\rho(m, m) \neq \rho(n, n)$. Now because $\Lambda_{BD}$ is strongly connected there exist $\beta \in s(\alpha) \Lambda_{BD} r(\rho)$. To finalise the proof we show that $\lambda = \alpha \beta \rho$ satisfies (\ref{APC}): Since
        \begin{align*}
                s \Big( \lambda \big( m, m + d(\lambda) - (m \vee n) \big) \Big) &= s \Big( (\alpha \beta \rho) \big( m, m + d(\alpha \beta \rho) - (m \vee n) \big) \Big) \\
                                &= s \Big( (\alpha \beta \rho) \big( m, m + d(\alpha \beta) \big) \Big) \\
                                &= \rho(m, m) \\
                                &\neq \rho(n, n) =  s \Big( \lambda \big( n, n + d(\lambda) - (m \vee n) \big) \Big)  
        \end{align*}
we see that $\lambda \big( m, m + d(\lambda) - (m \vee n) \big) \neq \lambda \big( n, n + d(\lambda) - (m \vee n) \big)$.
		
		If $BD$ instead has a red $a$-breaking cycle then a similar argument can be constructed by reversing the roles of the coordinates.
        \end{proof}

        \begin{remark*}
In Examples \ref{rem:b4APT} (1) we showed that the basic data $BD$ from Example~\ref{ldr} (2) has a blue $0$-breaking cycle and so
satisfies the hypotheses of Theorem~\ref{APT}. Hence we are able to prove aperiodicity for $2$-graphs which do not arise from
the parameters of \cite{prw}.
\end{remark*}

\noindent The aperiodicity result given in \cite[Theorem~5.2]{prw} employs a condition, called three invertible corners. It implies that all
$2$-graphs which are defined using tiles $T$ with $c_1 =0$ or $c_2=0$ cannot be aperiodic (see \cite{prw2}). Using Weaver's argument from \cite[Lemma 6.3.1]{W}, we may show that the same is true here:

\begin{lemma}
        Let $BD$ be basic data for a tile with $c_1=0$ or $c_2=0$, then $BD$ does not have an $a$-breaking cycle for all $a \in A$. Moreover,
        $\Lambda_{BD}$ is not aperiodic and so $C^* ( \Lambda_{BD} )$ is not simple.
        \end{lemma}

        \begin{proof}
        Suppose $c_2=0$ and $BD$ is basic data for the tile $T=\{0, e_1 , \ldots , c_1 e_1 \}$.
        Since $T \cap ( T + e_2 ) = \emptyset$, there is no blue $a$-breaking cycle for any $a \in A$.

        Suppose $v \in \Lambda_{BD}^0$ is such that $v |_{T \cap ( T + e_1 )} \equiv a$ for some $a \in A$.
        Since $T \cap ( T+e_1) = \{ e_1 , \ldots , c_1 e_1 \}$ this implies that $v(  j e_1 ) = a$ for all $1 \le j \le c_1$.
        By \eqref{eq:Fpdef} $v(0)=v(c_2 e_2)$ is completely determined by $v (c_1 e_1 )$, so there is only one vertex
        $v \in \Lambda_{BD}^0$ with $v |_{T \cap ( T + e_1 )} \equiv a$ for some $a \in A$.
        Hence by Definition \ref{ac} there is no red $a$-breaking cycle for any $a \in A$.


        Since $c_2=0$ it follows by Proposition \ref{sk} (3) that $| v \Lambda_{BD}^{e_1} | = 1$ for all $v \in \Lambda_{BD}^0$;
        in particular there is only one blue edge coming into each $v \in \Lambda_{BD}^0$.
        Since $\Lambda_{BD}^0$ is finite, it follows that every sufficiently long blue path in $\Lambda_{BD}$ must pass through the same vertex twice. Hence the blue subgraph of the skeleton of $\Lambda_{BD}$ consists of a finite number of disjoint cycles.
        Let $n$ be the lowest common multiple of the lengths of these cycles.
        Hence every path $\lambda \in \Lambda_{BD}$ has a horizontal repeating pattern of a block of width $n$. In particular for $\lambda \in \Lambda_{BD}$ and $l \in \Z^2$ we have $\lambda\addon_{T+l} = \lambda\addon_{T+l+ne_1}$ whenever $T+l$ and $T+l+ne_1$ are in the domain of
        $\lambda$. Fix $v \in \Lambda_{BD}^0$, $m,p \in \N^2$ satisfying $p_1-m_1=n$ and $p_2-m_2=0$. Then every $\lambda \in v \Lambda_{BD}$ with
        $d(\lambda)\ge p$ has
        \[
        \lambda ( m , m + d ( \lambda ) -p ) = \lambda ( p , d ( \lambda ) ) .
        \]

        \noindent Hence \eqref{APC} cannot be satisfied which implies that $\Lambda_{BD}$ is not aperiodic by Lemma \ref{lem:reducecase}. Thus \cite[Theorem 3.1]{robsi1} implies that $C^* ( \Lambda_{BD} )$ is not simple.

        Similar arguments apply for tiles with $c_1=0$.
        \end{proof}

\noindent For a general tile, we do not know if the converse of Theorem \ref{APT} holds. Recall from Proposition~\ref{bdreln} that given parameters $( T, q, 0 , w)$ from \cite{prw} with $w(c_1 e_1 )$ and $w ( c_2 e_2 )$
invertible in $\Z / q \Z$ then there is basic data $BD$ such that $\Lambda_{BD} \cong \Lambda ( T, q, 0 , w)$. We now show that if $w(0)$ is invertible in $\Z / q \Z$ (and so $\Lambda (T,q,0,w)$ is aperiodic by \cite[Theorem~5.2]{prw}) then $\Lambda_{BD}$ is aperiodic.

\begin{prop} \label{prop:subsume}
Suppose $(T,q,0,w)$ is a set of parameters and  $BD$ is a set of basic data such that $\Lambda_{BD} \cong \Lambda ( T, q , 0 , w)$. If $c_1 , c_2 \ge 1$ and $w(0)$ is invertible in $\Z / q \Z$ then $BD$ has an $a$-breaking cycle for some $a \in \Z / q \Z$ and so $\Lambda_{BD}$ aperiodic.
\end{prop}

\begin{proof}

We claim that the basic data $BD$ has a blue $1$-breaking cycle. We do this by finding distinct vertices $v_1 , \ldots , v_k \in \Lambda_{BD}^0$ with $v_i |_{T \cap ( T + e_2 )} \equiv 1$ which satisfy either (1) or (2) of Definition \ref{ac}.
Let $n = \sum_{i \in T \cap ( T + e_2 )} w(i)$ and $n' \in \Z / q \Z$ the unique number such that $n+n'=0$. Since vertices $v \in \Lambda (T,q,0,w)^0$
must satisfy $\sum_{i \in T} v (i ) w(i) = 0$, the vertices we
produce must therefore have a weighted sum of $n'$ along the bottom row.

If $n'=0$, then since $w(0), w(c_1 e_1)$ are invertible in $\Z / q \Z$ we may find nonzero elements $a,b \in \Z / q \Z$ such that $a w(0) + b w (c_1 e_1)=0$. Define  $v_1 : T \to \Z / q \Z$ by $v_1 |_{T \cap ( T+e_2)} \equiv 1$, with $v_1 ( 0 ) = a$, $v_1 ( c_1 e_1 ) = b$ and if $c_1 \ge 2$ we set $v_1 ( j e_1 ) = 0$ for $1 \le j \le c_1-1$. Since the weighted sum along the bottom row is $0$ we have $v_1 \in \Lambda (T,q,0,w)^0 = \Lambda_{BD}^0$.

If $n' \neq 0$. Define $v_1 : T \to \Z / q \Z$ by $v_1 |_{T \cap ( T+e_2)} \equiv 1$, with $v_1 ( j e_1 ) = 0$ for $1 \le j \le c_1$ and $v_1 ( 0 ) = n' w(0)^{-1}$. Since the weighted sum along the bottom row is $n'$ we have $v_1 \in \Lambda (T,q,0,w)^0
= \Lambda_{BD}^0$.

Define $v_2 : T \to \Z / q \Z$ by $v_2 |_{T \cap ( T+e_2)} \equiv 1$, with $v_2 ( j e_1 ) = v_1 ( (j+1 ) e_1 )$ for $0 \le j \le c_1-1$ and $v_2 ( c_1 e_1 ) = w ( c_1 e_1 )^{-1} \left( n' - \sum_{j=0}^{c_1-1} w ( j e_1 ) v_2 ( j e_1 ) \right)$. Since the weighted sum along the bottom row is $n'$ we have $v_2 \in \Lambda (T,q,0,w)^0 = \Lambda_{BD}^0$ and $| v_1 \Lambda^{e_1} v_2 | = 1$.

If $n' \neq  0$ then one checks that $v_1 \neq v_2$. If $n'=0$ then $v_1 = v_2$ only if $c_1 =1$ and $a=b$; in this case we define $v_2 : T \to \Z / q \Z$ by $v_2 |_{T \cap ( T+e_2)} \equiv 1$, with $v_2 ( j e_1 ) = 0$ for $0 \le j \le c_1$. Since the weighted sum along the bottom row is $0$ we have $v_2 \in \Lambda (T,q,0,w)^0 = \Lambda_{BD}^0$ and $| v_2 \Lambda^{e_1} v_2 | = 1$. Hence the vertices $v_1 , v_2$ satisfy condition (2) of Definition \ref{ac}. Hence we may assume that $v_1 \neq v_2$ from here on.

Now define $v_3 : T \to \Z / q \Z$ by $v_3 |_{T \cap ( T+e_2)} \equiv 1$ by $v_3 ( j e_1 ) = v_2 ( (j+1) e_1 )$ for $0 \le j \le c_1 -1$, and $v_3 ( c_1 e_1 ) = w ( c_1 e_1 )^{-1} \left( n' - \sum_{j=0}^{c_1-1} w ( j e_1 ) v_3 ( j e_1 ) \right)$. Since the weighted sum along the bottom row is $n'$ we have $v_3 \in \Lambda ( T,q,0,w)^0$ and $| v_2 \Lambda_{BD}^0 v_3 | = 1$.

If $v_3 = v_1$  then $\{ v_1 , v_2 \}$ satisfy condition (1) of Definition \ref{ac}. Suppose that $v_2 = v_3$ then entries of the bottom row of $v_2$ must all be equal. But $v_2 (0)=0$ and $v_2 ( c_1 e_1 ) = w ( c_1 e_1 )^{-1} \left( n' - \sum_{j=0}^{c_1-1} w ( j e_1 ) v_2 ( j e_1 ) \right) = w ( c_1 e_1 )^{-1} n' \neq 0$ which is a contradiction.
Hence $v_2 \neq v_3$.

Continuing in this way, we produce $k \ge 2$ vertices $v_1 , v_2 , \ldots , v_k$ in $\Lambda (T,q,0,w) = \Lambda_{BD}^0$ with $| v_i \Lambda v_{i+1}| = 1$ for $1 \le i \le k-1$. Since $\Lambda_{BD}^0$ is finite we must have $v_i = v_j$ for some $i < j$. Hence the vertices $v_i , \ldots , v_j$ satisfy condition (1) of Definition \ref{ac}. The final statement now follows from Theorem \ref{APT}.
\end{proof}

\noindent The proof of Proposition \ref{prop:subsume} may be modified to show that there are also red and blue $p$-breaking cycles for all $p \in \Z / q \Z$.

\begin{remark*}
It follows from Proposition \ref{prop:subsume} that our aperiodicity result, Theorem \ref{APT}, subsumes  \cite[Theorem 5.2]{prw}. Furthermore, Theorem~\ref{APT} allows us to deduce aperiodicity for a wider class of $2$-graphs described using the parameters of \cite{prw}. In Examples \ref{rem:b4APT} (3) we gave basic data $BD$ which has a blue $0$-breaking cycle, and so by Theorem \ref{APT} we conclude that the associated $2$-graph is aperiodic. However the basic data $BD$ corresponds to the parameters $( T,2,0,w )$ of \cite{prw} where the rule $w$ is given by $w(0)=0$ and $w(e_1)=w(e_2)=w(2e_1)=1$. Since $w(0)$ is not invertible in $\Z / 2 \Z$, we cannot apply \cite[Theorem~5.2]{prw} to deduce aperiodicity.
\end{remark*}


        \section{Simplicity of basic data $2$-graph $C^*$-algebras.}\label{simplicity}

        \noindent We can now use the properties of our basic data $2$-graphs to determine the structure of the corresponding higher rank graph $C^*$-algebras. Recall from \cite[Remark~A.3]{ls} that for a row-finite 2-graph $\Lambda$ with no sources, we say $\Lambda$ is \emph{cofinal} if for every $v, w \in \Lambda^0$ there exist $n \in \N^2$ such that $v \Lambda s(\alpha) \neq \emptyset$ for every $\alpha \in w\Lambda^n$.

        \begin{theorem}\label{class}

       Let $\Lambda_{BD}$ be the 2-graph associated to the basic data $BD$ which has an $a$-breaking cycle. Then $C^*(\Lambda_{BD} )$ is unital, nuclear, simple, purely infinite and belongs to the bootstrap class $\mathcal{N}$.

        \end{theorem}
        \begin{proof}

        Since $\Lambda_{BD}^0$ is finite the sum $\sum_{v \in \Lambda_{BD}^0} s_v$ is an identity for $C^* (\Lambda_{BD} )$. As $\Lambda_{BD}$ is a 2-graph, \cite[Theorem~5.5]{kp} implies $C^* (\Lambda_{BD} )$ belongs to the bootstrap class $\mathcal{N}$ for which the UCT holds, and so $C^*(\Lambda_{BD} )$ is nuclear. Since $\Lambda_{BD}$ has an $a$-breaking cycle, $\Lambda_{BD}$ is aperiodic by Theorem \ref{APT}, and hence \cite[Lemma 3.2]{robsi1} implies $\Lambda_{BD}$ has no local periodicity. Since $\Lambda_{BD}$ is strongly connected it is clearly cofinal. Then \cite[Remark~A.3 and Theorem 5.1]{ls} show that the definition of cofinality used here is equivalent to the definition of cofinality used in \cite{robsi1}. Then \cite[Theorem 3.1]{robsi1} implies that $C^*(\Lambda_{BD} )$ is simple.

        We use the results of \cite{sims2} to show that $C^*(\Lambda_{BD} )$ is purely infinite: we must show that every $v \in \Lambda_{BD}^0$ can be reached from a loop with an entrance and that $\Lambda_{BD}$ satisfies the technical condition (C) of \cite[Section 7]{sims2}.  The argument that every $v \in \Lambda_{BD}^0$ can be reached from a loop with an entrance is the same as in the proof of \cite[Theorem~6.1]{prw}. It remains to show that $\Lambda_{BD}$ satisfies (1) and (2) of condition (C). Since $\Lambda_{BD}$ is aperiodic we may use \cite[Proposition 3.6]{ls} to conclude that $\Lambda_{BD}$ satisfies (1) of condition (C). By Proposition \ref{sk}, $\Lambda_{BD}$ is row-finite with no sinks or sources, so each finite exhaustive set is a union of sets of the form $\Lambda_{BD}^n$ and hence trivially satisfies (2) of condition (C).  The result now follows from \cite[Proposition~8.8]{sims2}.
        \end{proof}

 \section{Multi-dimensional shift spaces of basic data 2-graphs.}\label{ss}

        Here we connect our work to the literature on dynamical systems. In our construction of the $2$-graph $\Lambda_{BD}$ from the basic data $BD$ we have have not used any algebraic properties of $A$ in our definitions. For this reason we align our work to that of Quas and Trow~\cite{qt} rather than the work of Schmidt~\cite{schm}. Given basic data $BD$ and associated $2$-graph $\Lambda_{BD}$, the two-sided infinite path
        space $\Lambda_{BD}^\Delta$ may be given the structure of a shift space (see \cite{kp2}). In Theorem~\ref{hmeo} we show that the shift space  to associated 2-graph $\Lambda_{BD}$ is homeomorphic to a shift of finite type $Y_{BD}$, as described in \cite{qt}.

        Fix a set of basic data $BD$ and define $Y_{BD}$ to be the subset
        \[
                \{ y \in A^{\Z^2} : y\addon_{T + n} \in \Lambda_{BD}^0 ~ \text{for every} ~ n \in \Z^2 \} \subset A^{\Z^2}.
        \]
Then it is routine to show that $Y_{BD}$ is a shift space with shift map $\sigma_b$  given by $( \sigma_b y ) (i) = y(i+b)$for $b \in \Z^2$ . In fact $Y_{BD}$ is a shift of finite type since $Q = \Lambda_{BD}^0$ is a set of allowed finite configurations for $S = T \subset \Z^2$.

        Kumjian and Pask~\cite{kp2} introduced the idea of a two-sided infinite path space for a $k$-graph $\Lambda$ with no sinks or sources. We will use their construction with $k = 2$. Let
        \[
                \Delta = \{ (m, n) : m, n \in \Z^2, m \leq n \};
        \]
where $r(m, n) = m$, $s(m, n) = n$ and $d(m, n) = n - m$. Then $(\Delta, d)$ is a 2-graph with no sources or sinks. Then the two-sided infinite path space of $\Lambda$ is
        \[
                \Lambda^\Delta = \{ x : \Delta \rightarrow \Lambda ~ \text{such that} ~ x ~ \text{is a degree preserving functor} \}.
        \]
Kumjian and Pask show that $\Lambda^\Delta$ has a locally compact topology generated by the cylinder sets
        \[
                Z(\lambda, n) = \{ x \in \Lambda^\Delta : x \big( n, n + d(\lambda) \big) = \lambda \} ~ \text{where $\lambda \in \Lambda$ and $n \in \Z^2$.}
        \]
In \cite[Section 3]{kp2} Kumjian and Pask also demonstrated that there is a metric $\rho$ which induces this topology and that $\Lambda^\Delta$ is compact if $\Lambda^0$ is finite. This fact is used in the proof of Theorem~\ref{hmeo} and details of this construction are given in \cite[Section~3]{kp2} and \cite[Section~1.1]{sz}. For every $b \in \Z^2$, define $\sigma^b : \Lambda^\Delta \rightarrow \Lambda^\Delta$ by $\sigma^b (x) (m, n) = x(m + b, n + b)$. A straightforward argument shows that $\sigma^b$ is a homeomorphism.

        \begin{theorem}\label{hmeo}

        Suppose we have basic data $BD$ with associated $2$-graph $\Lambda_{BD}$, and $Y_{BD}$ defined as above. Then there is a homeomorphism $h : \Lambda_{BD}^\Delta \rightarrow Y_{BD}$ such that $\sigma_b \circ h = h \circ \sigma^b$ for every $b \in \Z^2$.

        \end{theorem}
        \begin{proof}

        Recall for any $x \in \Lambda_{BD}^\Delta$ that $x(i, i) \in \Lambda_{BD}^0$ is a function from $T$ to $A$, so in order to evaluate $x(i)$ we evaluate $x(i, i)$ at 0. Define $h : \Lambda_{BD}^\Delta \rightarrow A^{\Z^2}$ by $\big( h(x) \big) (i) = x(i, i)(0)$ where $x \in \Lambda^\Delta$ and $i \in \Z^2$. Fix $x \in \Lambda_{BD}^\Delta, n \in \Z^2$ and let $p = x|_{P + n}$. To show $h : \Lambda_{BD}^\Delta \rightarrow Y_{BD}$ we must show $h(x)\addon_{T + n} \in \Lambda_{BD}^0$. Since $p \in A^P$ it only remains to show $\big( h(x) \big (n + c_1 e_1) = f_p \big( x(n + c_2 e_2) \big)$. Since $x(n, n) \in \Lambda_{BD}^0$ we have
        \[
                \big( h(x) \big (n + c_1 e_1) = x(n + c_1 e_1, n + c_1 e_1)(0) = x(n, n)(c_1 e_1) = f_p \big( x(n + c_2 e_2) \big).
        \]

        Next we show $h$ is a homeomorphism. Since $\Lambda_{BD}^0$ is finite it follows that $\Lambda_{BD}^\Delta$ is compact and so it suffices show that $h$ is a continuous bijection. To this end for each $y \in Y_{BD}$ we define a functor $g(y) \in \Lambda_{BD}^\Delta$ by
        \[
                \big( g(y) \big) (m, n) = y|_{T(n - m) + m} ~ \text{for} ~m, n \in \N^2 ~ \text{such that} ~ m \leq n.
        \]
Then the definition of $Y_{BD}$ ensures $\big( g(y) \big) (m, n)$ is a path in $\Lambda_{BD}$ with degree $n - m$. So $g$ is degree preserving. To see that $g(y)$ is a functor from $\Delta$ to $\Lambda_{BD}$, fix $m, M \in \Z^2$ such that $m \leq n \leq M$. Then Proposition~\ref{comp} implies
        \begin{align*}
                \big( g(y) \big) (m, M) = y|_{T(M - m) + m} = y|_{T(n - m) + m} y|_{T(M - n) + n} = \big( g(y) \big) (m, n) \big( g(y) \big) (n, M).
        \end{align*}
Since $g(y)$ is a degree preserving functor for every $y \in Y_{BD}$ we now know that $g : Y_{BD} \rightarrow \Lambda_{BD}^\Delta$ and we claim $g = h^{-1}$. Fix $y \in Y_{BD}$ and $i \in \N^2$. Then
        \[
                \Big( h \big( g(y) \big) \Big) (i) = \big( g(y) \big) (i, i) (0) = y|_{T(0) + i} (0) = y(i)
        \]
To finalise the claim we must show $g \big( h(x) \big) = x$. Fix $m, n \in \N^2$ such that $m \leq n$ and let $i \in T(n - m)$. Then
        \begin{align*}
                \Big( \big( g(h(x)) \big) (m, n) \Big) (i) &= \big( h(x)|_{T(n - m) + m} \big) (i) \\
                                &= h(x) (i + m) \\
                                &= x(i + m, i + m)(0)
                                = x(i + m, i + n)(0)
                                = x(m, n)(i).
        \end{align*}
This proves the claim that $g = h^{-1}$. Thus to show that $h$ is a homeomorphism it remains to show that $h$ is continuous. Since $\Lambda_{BD}^\Delta$ and $Y_{BD}$ are both metric spaces it suffices to show that $h(x_m)(i) \to h(x)(i)$ for every $i \in \Z^2$. Fix $i \in \Z^2$. Then since $\sigma^{-i} (x) = x$ it follows that $Z \big( x(i, i), 0 \big)$ is an open neighbourhood of $x$. So for sufficiently large $m$ we have $x_m \in Z \big( x(i, i), 0 \big)$. That is, $\lim_{m \rightarrow \infty} x_m (0, 0) = x(i, i)$ and hence
        \[
                \lim_{m \rightarrow \infty} h(x_m)(i) = \lim_{m \rightarrow \infty} x_m (i, i)(0) = x(i, i)(0) = h(x)(i).
        \]

        To finalise the proof for $x \in \Lambda_{BD}^\Delta$ and $b \in \Z^2$ we compute
        \[
                \sigma_b \big( h(x) \big) (i) = \big( h(x) \big) (i + b) = x(i + b, i + b)(0) = \sigma^b (x) (i, i)(0) = h \big( \sigma^b (x) \big) (i). \qedhere
        \]
        \end{proof}

        \begin{remark} \label{rem:T03} If $T=\{0\}$ then the $2$-graph associated to the basic data $BD=(T,A,\{f_a\})$ has only one infinite path, the function which is identically $a$. One checks that this function is also the only element of $Y_{BD}$.
        \end{remark}

        \noindent Given a two dimensional shift space $X$, let $\mathcal{B}_d(X)$ denote the set of all $d \times d$ configurations which occur in $X$. As in \cite{qt} we define the \emph{topological entropy} of a shift space $X$ to be $\lim_{d \rightarrow \infty} \frac{1}{2^d} \log |\mathcal{B}_d(X)|$ and denote it by $h(X)$. Suppose we have a set of basic data $BD$. We now compute the topological entropy of $\Lambda_{BD}$ (and hence $Y_{BD}$).

By definition $|\mathcal{B}_d( \Lambda_{BD} )| = |\Lambda_{BD}^{(d, d)}|$ for all $d \ge 1$. Since $|\mathcal{B}_d( \Lambda_{BD} )| \leq |A|^{d^2}$ for all $d \ge 1$ the topological entropy of $\Lambda_{BD}$ is,
        \begin{align*}
                h( \Lambda_{BD}) = \lim_{d \rightarrow \infty} \frac{1}{2^d} \log |\mathcal{B}_d( \Lambda_{BD} )| \leq \lim_{d \rightarrow \infty} \frac{1}{2^d} \log |A|^{d^2} = \lim_{d \rightarrow \infty} \frac{d^2}{2^d} \log |A| = 0.
        \end{align*}
This agrees with \cite[Proposition~1.2]{sz} where a different definition of entropy, formulated using the notion of separating subsets, is used.

    \end{document}